%% file: article.tex
\documentclass[11pt,a4paper]{article}
\usepackage[symbol]{footmisc}
 
\usepackage{amsmath, amsfonts, amssymb, amsbsy}
\usepackage{latexsym} \usepackage{stmaryrd} 
\usepackage{accents}

\usepackage[textsize=footnotesize]{todonotes}

\usepackage[labelfont=bf]{caption}
\usepackage{hyperref}
\hypersetup{
    colorlinks,%
    citecolor=black,%
    filecolor=black,%
    linkcolor=black,%
    urlcolor=black
}

\newcommand{\bl}{\color{black}}

\begin{document}

\input{macro} 

\input{cover}

\section*{Introduction}

Griffith's criterion \cite{Griffith20} was originally proposed for the quasi-static propagation of straight cracks in brittle materials. Mathematically, it is a rate-independent system comprising (a) an irreversibility constraint, (b) an equilibrium condition, and (c) an evolution law (or flow rule). Usually it is expressed by means of Karush-Kuhn-Tucker complementarity conditions, in terms of toughness and \textit{energy release} or, equivalently, in terms of derivatives of the energy, or stress intensity factors. 

In order to briefly describe the model, we 
assume for the moment that $\Omega \subset \mathbb{R}^2$ is the reference configuration and consider a family of straight cracks $K_\ell \subset \Omega$ parametrized by their length $\ell$. Let 
$$	
	E (u, K_\ell) = \int_{\Omega \setminus K_\ell}  W (\strain (u)) \, dx	
$$
be the elastic energy for a displacement $u \in H^1( \Omega \setminus K_\ell ; \mathbb{R}^2)$. 
Given a time depending boundary datum and a crack $K_\ell$, the {\it elastic energy} at equilibrium is given by 
\begin{equation} \label{reduced} 
	\tilde{E} (t, \ell) = \min \{ E (u, K_\ell) : u \in H^1 ( \Omega \setminus K_\ell; \mathbb{R}^2) \text{ with }  u = g(t) \text{ on $\partial_D \Omega$} \}  .   
\end{equation}
The {\it dissipated energy} is instead proportional to crack elongation, and takes the form $ D (\ell ) = G_c (\ell - \ell_0)$,
where $G_c >0$ is the toughness while $\ell_0$ is the length of the initial crack. The {\it total energy} of the system is then given by 
$$ F (t , \ell) = \tilde{E} (t, \ell) + \mathcal{D} ( \ell) .   $$
The {\it energy release} is the (configurational) variation of the elastic energy $\tilde{E}$ with respect to crack growth, which in this setting boils down to the derivative $$ G (t, \ell) = - \partial_\ell \tilde{E} ( t, \ell) . $$

\medskip
With this notation, the Karush-Kuhn-Tucker conditions, which give the evolution in classical form, are
\begin{equation} \label{KKT} 	
	G ( t, \ell(t) ) \le G_c  , \quad \big( G ( t, \ell(t) ) - G_c \big) \, \dot{\ell} (t) = 0 ,
\end{equation} 
together with the monotonicity constraint $\dot{\ell}(t) \ge 0$, which models irreversibility. 
In terms of the energy $F$, they read
\begin{equation} \label{KKT-F}
	\partial_\ell F (t, \ell(t)) \ge 0  , \quad  \partial_\ell F (t, \ell(t)) \, \dot{\ell} (t) = 0 .	
\end{equation} 
A few comments are due. The first term gives the equilibrium with respect to the admissible variations, i.e.~for growing cracks, by irreversibility. The latter equation governs the evolution, asserting the following:
\begin{itemize} 
\vspace{-5pt}\item the crack does not grow if the energy release is strictly lower than the toughness,
\vspace{-5pt}\item the crack grows when the energy release equals the toughness.
\end{itemize}

At this point, it is important to notice that the KTT system above holds only for steady state (or stable) evolutions, while it breaks when the energy release gets strictly greater than the toughness. 
In this case, unsteady (instantaneous) propagation occurs. Mathematically, the evolution is instantaneous and some sort of ``extended criterion'' is needed. Different solutions have been proposed in the literature. 

Following \cite{NegriOrtner_M3AS08}, in the unstable times $t$ we may assume that $G ( t , l) \ge G_c$ for every $l \in [ \ell_-(t) , \ell_+(t) ] $. 
In other terms, the time $t$ is a discontinuity point where the evolution jumps from $\ell_-(t)$ to $\ell_+(t)$ passing only through unstable configurations; typically, it happens that $G (t, \ell_- (t)) = G_c$,  $G (t, l ) > G_c$ for $l \in ( \ell_-(t) , \ell_+(t) )$, and $G (t, \ell_+ (t)) = G_c$. 

Alternatively, we may employ {\it energetic solutions} \cite{FrancMar98,Giac05,MielkeRoubicek}. These evolutions, based on unilateral global minimizers of the energy, satisfy the KKT system in continuity points; they may as well present discontinuities, however, their qualitative behaviour is different from that of \cite{NegriOrtner_M3AS08}. Indeed, in the discontinuity times 
the transition between $\ell_-(t)$ and $\ell_+(t)$ is not fully unstable. In other terms, it may happen (for instance at crack initiation) that $G ( t, \ell_- (t) ) < G_c$ and that $G (t , l ) < G_c$ for some $l \in [\ell_- (t) , \ell_+(t)]$, while $G ( t, \ell_+ (t) ) = G_c$. On the contrary, the KKT conditions \eqref{KKT} are satisfied in continuity points. 

Let us mention that the unstable regimes could be described out of the quasi-static framework using rate-dependent and dynamical models, see e.g.~\cite{LarsenOrtnerSuli_M3AS10, Roub20, FTSWW}.

\medskip
When the path of the crack is not known a priori, the system of KKT conditions is complemented by a {\it directional criterion}, whose role is to select the "geometry" of the crack increment. Several criteria exist in the literature, for instance, the {\it principle of local symmetry} and the {\it maximum energy release rate}, which give different but comparable predictions. 
Here, we mention only the latter, which is closely related to the phase-field framework. For sake of simplicity, we do not enter into the fine technical details (interesting but out of scope here) and we confine ourselves to steady-state evolutions. Let $K \subset \Omega$ be a crack and consider a system of coordinates centred in its tip. For $\zeta \in \mathbb{R}^2$ with $| \zeta |=1$, the variations of $K$ in direction $\zeta$ take the form $K \cup [0, h \zeta]$ where $[0,h\zeta]$ is the line segment $\{ s \, \zeta : s \in [0, h] \}$. Then, the variation of energy in direction $\zeta$ is 
\begin{equation} \label{GMAX1}
	d E ( t, K ; \zeta ) = \lim_{h \to 0^+} \frac{\tilde{E} ( t, K \cup [0,h\zeta]  ) - \tilde{E} (t, K) }{h} ,
\end{equation}
where $\tilde{E} (t,K)$ denotes the elastic energy at equilibrium, cfr.~\eqref{reduced}. With this notation, the {\it maximal energy release rate} is given by 
\begin{equation} \label{GMAX2}
	G_{\text{\sl max} }\, ( t ,K) = \max \,\{ - d \tilde{E} (t, K; \zeta) : | \zeta | =1 \} , 
\end{equation} 
i.e the derivative of the elastic energy along the direction of steepest descent. 
Let $\gamma$ be a curve parametrizing in time the evolution of the tip. The crack set is then $K(t) = K_0 \cup \gamma ([0,t])$, where $K_0$ is the initial crack. Denoting by $\ell(t)$ the length of the crack $K(t)$ we have $\dot{\ell}(t) = | \dot\gamma(t)|$. In this framework, the set of (scalar) KKT conditions becomes
\begin{equation} \label{e.GR1}
	G_{\text{\sl max}}\, ( t , K(t) ) \le G_c , \quad \big( G_{\text{\sl max}}\, ( t , K(t) ) - G_c \big) \, \dot{\ell}(t) = 0, 
\end{equation} 
together with 
\begin{equation} \label{e.GR2}
d \tilde{E} (t, K (t); \dot \gamma (t) ) = - | \dot \gamma (t) | \, G_{\text{\sl max}} (t, K(t)) ,
\end{equation}  which selects the propagation direction as that of steepest descent for the elastic energy.

\medskip
Now, let us turn to the {\it phase-field} framework. In this setting, proposed in \cite{BourdFrancMar00}, the crack is represented by the phase-field variable $v \in [0,1]$, with $v=0$ and $v=1$ corresponding to fracture and sound material, respectively. The total energy takes the form 
\begin{equation}\label{F}
	\F_\eps  (u, v) = \E_\eps  (u,v) + G_c \, \mathcal{L}_\eps (v) , 
\end{equation} 
where $u$ is the displacement and $\eps>0$ is the internal length. 
Phase-field {\it evolutions} are computed using time discrete incremental problems:
let  $t^n_k = k T / n $ for $k=0,...,n$, given $v^n_{0}$, $( u^n_k, v^n_k)$ are such that: 
\begin{equation*} 
\begin{cases}
	u^n_{k} \in \argmin \{  \F_\eps  (  u , v^n_{k} ) \} & \text{for $k=0,...,n$,} 
\\
         v^n_{k} \in \argmin \{  \F_\eps ( u^n_{k} , v ) : v \le v^n_{k-1} \} &  \text{for $k=1,...,n$} . 
\end{cases} 
\end{equation*} 

Here, we do not specify the (sub)-scheme which drives $( u^n_{k-1}, v^n_{k-1} )$ to $( u^n_k, v^n_k )$. Note that, being the energy non-convex, in general the update $( u^n_k, v^n_k )$ is non-unique and may depend on the adopted (sub)-scheme. This framework is intentionally large, in order to include all the algorithms employed in the numerical simulations, such as staggered and monolithic \cite{BourdFrancMar00,FarrellMaurini_IJNME17}, and even energetic evolutions.

The first goal of our paper is the characterization of the time continuous limit. To this end we employ a linear interpolation of the configurations $(u^n_k, v^n_k)$ together with a suitable Lipschitz reparametrization in a finite time interval $I$. 
In the limit, stable and unstable branches of the evolution correspond to disjoint subsets $I_s$ and $I_u$ of $I$, characterized by Kuratowski convergence. First, we provide the following characterization, in terms of derivatives of the energy. For a.e.~$t \in I_s $ the configuration $(u(t), v(t))$ satisfies 
 \begin{align} \label{KKTphase}
	\partial_u  \F_\eps ( u (t) , v(t) ) [ \phi ] = 0, \quad \partial_v \F_\eps ( u (t) , v (t) ) [ \xi ] \ge 0 , \quad \partial_v \F_\eps ( u(t) ,  v (t) ) [ \dot{v}(t)] = 0 ,
\end{align}
for every admissible variation $\phi$ and $\xi$, which plays the role of \eqref{KKT-F}. Technically, in the proof of the above conditions there are a couple of delicate points: the strong convergence of the phase-field functions (see Lemma \ref{ls.strong}) and the proof of the power identity (see Theorem \ref{t.enbal}). 

Our next goal is the characterization of the limit evolution in terms of Griffith's criterion. In analogy with the sharp-crack setting (\ref{GMAX1}-\ref{GMAX2}) we first introduce the {\it energy release rate} as the variation of elastic energy with respect to the variations of crack length, i.e.
\begin{align*}
	\mathcal{G}_{\text{\sl max}}  (t,v) =  \biggl\{ \! - \!\lim_{\,\,h\to 0^+} \frac{\tilde{\E}_\eps  (t,v+h\xi)-\tilde{\E}_\eps (t,v)  }{\mathcal{L}_\eps (v+h\xi)-\mathcal{L}_\eps(v)}  \,:\,\xi \le 0  \biggr\}  , 
\end{align*} 
where $ \tilde{\E} ( t ,v)$ is again the elastic energy at equilibrium, and $\xi \le 0$ gives the admissible variations, by irreversibility. Denoting $\dot{\ell} ( t ) =d\mathcal{L}_\eps (v(t))[\dot{v}(t)]$ we have the following properties, providing Griffith's criterion for phase-field fracture: for steady state propagation, i.e.~for $t \in I_s$, we have 
$$
	\mathcal{G}_{\text{\sl max}} ( t  , v(t) ) \le G_c , \quad  \big( \mathcal{G}_{\text{\sl max}} ( t , v(t) ) - G_c \big) \, \dot{\ell} (t) = 0 
$$
together with steepest descent condition 
$$ d \tilde{\E}_\eps ( t, v(t)) [ \dot{v}(t) ] = - | \dot \ell (t) |\,   \mathcal{G}_{\text{\sl max}} ( t  , v(t) )  . $$
In the unsteady regime, we have in general that $\mathcal{G}_{\text{\sl max}}  (t, v(t)) \ge G_c$ on a subset of positive measure of $I_u$. Finally, it holds $\dot\ell(t)\ge 0$ in the whole parametrization interval $I$.

Last, but not least, let us mention that the phase-field approach relies on the rigorous base of {\it $\Gamma$-convergence} \cite{AmbrosioTortorelli_CPAM90,DalMaso93,Braid98} and specifically on the converge of the energies as the internal lenght $\eps$ vanishes. In our setting, (under suitable technical assumptions \cite{ChambolleContiFrancfort_ARMA18}) the phase-field energies $\F$ converge to the sharp crack energy 
$$
	\int_{\Omega \setminus S_u}  \mu   | \strain_d |^2 +    \kappa   | \strain_v^+ |^2  \, dx  + G_c \mathcal{H}^{1} (S_u)   , 
$$
where the displacement $u \in SBD (\Omega)$ satisfies the non-interpenetration condition $\llbracket u \rrbracket \cdot \nu \ge 0 $ on the discontinuity set $S_u$. It is well known that $\Gamma$-convergence provides the convergence of global minimizers, while not much is known on the convergence of equilibrium points, energy release rate, and evolutions. In this respect we mention \cite{N_AML13,SicMar13,BabadjianMillotRodiac_22,Maggiorelli24}. The convergence of phase-field to sharp-crack evolutions is proven in \cite{N_AML13,SicMar13}, at the price of restrictive assumptions on the geometry of the crack. The convergence of critical points of the Ambrosio-Tortorelli energy is contained in \cite{BabadjianMillotRodiac_22}, while  \cite{Maggiorelli24} presents an accurate numerical study (in the steady and unsteady propagation regimes) on the convergence of energies, energy release rate and evolutions, endorsing the theoretical results presented in this paper.

\tableofcontents


\input{notation}

\input{preliminaries}

\input{conv-time}

\input{griffith}

\input{pde}

\input{alternate}

\appendix
\input{appendix}

\noindent {\bf Acknowledgement.} 
The author acknowledges support from PRIN 2022 (Project no. 2022J4FYNJ), funded by MUR, Italy, and the European Union -- Next Generation EU, Mission 4 Component 1 CUP F53D23002760006.
The authors are members of GNAMPA - INdAM and M. Negri is research associate at IMATI - CNR. 

\bibliographystyle{plain}
\bibliography{ref} 

\end{document}

%% file: macro.tex
\newtheorem{theorem}{Theorem}[section]
\newtheorem{proposition}[theorem]{Proposition}
\newtheorem{corollary}[theorem]{Corollary}
\newtheorem{lemma}[theorem]{Lemma}
\newtheorem{definition}[theorem]{Definition}
\newtheorem{remark}[theorem]{Remark}
\newtheorem{example}[theorem]{Example}
\newtheorem{conjecture}[theorem]{Conjecture}

\newcommand{\proof}{\noindent \textbf{Proof. }}
\newcommand{\qed}{ \hfill {\vrule width 6pt height 6pt depth 0pt} \medskip }

\newcommand{\separe}{\medskip \centerline{\tt -------------------------------------------- } \medskip}
\newcommand{\passo}[1]{\separe \noindent{\tt #1} \medskip}

\newcommand{\nota}[1]{{\small \tt [#1]}}
\newcommand{\quadro}[1]{\medskip \noindent \framebox{\begin{minipage}[c]{\textwidth} #1 \end{minipage}}\medskip}

\renewcommand{\epsilon}{\varepsilon} \newcommand{\eps}{\varepsilon} 
\renewcommand{\det}{\mathrm{det}} \newcommand{\0}{ {\mbox{\tiny 0}} }
\newcommand{\stress}{\boldsymbol{\sigma}} \newcommand{\strain}{\boldsymbol{\epsilon}}
\newcommand{\argmin}{ \mathrm{argmin} \,}\newcommand{\argmax}{ \mathrm{argmax} \,}
\newcommand{\I}{ {\mbox{\tiny \rm I}} } \newcommand{\II}{ {\mbox{\tiny \rm II}} }
\def\Xint#1{\mathchoice
 {\XXint\displaystyle\textstyle{#1}}%
 {\XXint\textstyle\scriptstyle{#1}}%
 {\XXint\scriptstyle\scriptscriptstyle{#1}}%
 {\XXint\scriptscriptstyle\scriptscriptstyle{#1}}%
 \!\int}
\def\XXint#1#2#3{{\setbox0=\hbox{$#1{#2#3}{\int}$}
 \vcenter{\hbox{$#2#3$}}\kern-.5\wd0}}
 \def\ddashint{\Xint=} 
 \def\dashint{\Xint-}  
\def\interior{\mathaccent'27}

\newcommand{\NOTE}[1]{\todo[inline, color=yellow!20]{\tt #1}}



\newcommand\ubar[1]{\underaccent{\bar}{\raisebox{-0.4pt}{\raisebox{0.4pt}{$#1$}}}}

\def\utilde#1{\mathord{\vtop{\ialign{##\crcr
$\hfil\displaystyle{#1}\hfil$\crcr\noalign{\kern1.5pt\nointerlineskip}
$\hfil\tilde{}\hfil$\crcr\noalign{\kern1.5pt}}}}}





\newcommand{\U}{\mathcal{U}} 
\newcommand{\Z}{\mathcal{Z}}
\newcommand{\E}{\mathcal{E}}
\newcommand{\F}{\mathcal{F}}
\newcommand{\D}{\mathcal{D}}
\renewcommand{\P}{\mathcal{P}}
\newcommand{\weakto}{\rightharpoonup}\newcommand{\weakstarto}{\stackrel{*}{\rightharpoonup}}
\newcommand{\R}{\mathbb{R}}
\newcommand{\V}{\mathcal{V}}

\newcommand{\wt}{\widetilde}
\def\calA{{\mathcal A}} \def\calB{{\mathcal B}} \def\calC{{\mathcal C}}
\def\calD{{\mathcal D}} \def\calE{{\mathcal E}} \def\calF{{\mathcal F}}
\def\calG{{\mathcal G}} \def\calH{{\mathcal H}} \def\calI{{\mathcal I}}
\def\calJ{{\mathcal J}} \def\calK{{\mathcal K}} \def\calL{{\mathcal L}}
\def\calM{{\mathcal M}} \def\calN{{\mathcal N}} \def\calO{{\mathcal O}}
\def\calP{{\mathcal P}} \def\calQ{{\mathcal Q}} \def\calR{{\mathcal R}}
\def\calS{{\mathcal S}} \def\calT{{\mathcal T}} \def\calU{{\mathcal U}}
\def\calV{{\mathcal V}} \def\calW{{\mathcal W}} \def\calX{{\mathcal X}}
\def\calY{{\mathcal Y}} \def\calZ{{\mathcal Z}}
\def\bbA{{\mathbb A}} \def\bbB{{\mathbb B}} \def\bbC{{\mathbb C}}
\def\bbD{{\mathbb D}} \def\bbE{{\mathbb E}} \def\bbF{{\mathbb F}}
\def\bbG{{\mathbb G}} \def\bbH{{\mathbb H}} \def\bbI{{\mathbb I}}
\def\bbJ{{\mathbb J}} \def\bbK{{\mathbb K}} \def\bbL{{\mathbb L}}
\def\bbM{{\mathbb M}} \def\bbN{{\mathbb N}} \def\bbO{{\mathbb O}}
\def\bbP{{\mathbb P}} \def\bbQ{{\mathbb Q}} \def\bbR{{\mathbb R}}
\def\bbS{{\mathbb S}} \def\bbT{{\mathbb T}} \def\bbU{{\mathbb U}}
\def\bbV{{\mathbb V}} \def\bbW{{\mathbb W}} \def\bbX{{\mathbb X}}
\def\bbY{{\mathbb Y}} \def\bbZ{{\mathbb Z}}
\def\FG{\mathbf} 
\def\mdot{\FG{:}}  \def\vdot{\FG{\cdot}} 
\def\bfa{{\FG a}} \def\bfb{{\FG b}} \def\bfc{{\FG c}} 
\def\bfd{{\FG d}} \def\bfe{{\FG e}} \def\bff{{\FG f}} 
\def\bfg{{\FG g}} \def\bfh{{\FG h}} \def\bfi{{\FG i}}
\def\bfj{{\FG j}} \def\bfk{{\FG k}} \def\bfl{{\FG l}} 
\def\bfm{{\FG m}} \def\bfn{{\FG n}} \def\bfo{{\FG o}} 
\def\bfp{{\FG p}} \def\bfq{{\FG q}} \def\bfr{{\FG r}} 
\def\bfs{{\FG s}} \def\bft{{\FG t}} \def\bfu{{\FG u}} 
\def\bfv{{\FG v}} \def\bfw{{\FG w}} \def\bfx{{\FG x}} 
\def\bfy{{\FG y}} \def\bfz{{\FG z}} 
\def\bfA{{\FG A}} \def\bfB{{\FG B}} \def\bfC{{\FG C}} 
\def\bfD{{\FG D}} \def\bfE{{\FG E}} \def\bfF{{\FG F}} 
\def\bfG{{\FG G}} \def\bfH{{\FG H}} \def\bfI{{\FG I}} 
\def\bfJ{{\FG J}} \def\bfK{{\FG K}} \def\bfL{{\FG L}} 
\def\bfM{{\FG M}} \def\bfN{{\FG N}} \def\bfO{{\FG O}} 
\def\bfP{{\FG P}} \def\bfQ{{\FG Q}} \def\bfR{{\FG R}} 
\def\bfS{{\FG S}} \def\bfT{{\FG T}} \def\bfU{{\FG U}} 
\def\bfV{{\FG V}} \def\bfW{{\FG W}} \def\bfX{{\FG X}} 
\def\bfY{{\FG Y}} \def\bfZ{{\FG Z}}

\newcommand{\abs}[1]{\left\lvert#1\right\rvert}      
\newcommand{\babs}[1]{\bigl\vert#1\bigr\vert} 
\newcommand{\Abs}[1]{\Bigl\vert#1\Bigr\vert} 
\newcommand{\Babs}[1]{\Bigl\vert#1\Bigr\vert} 
 \newcommand{\norm}[1]{\left\lVert#1\right\rVert}
\newcommand{\bnorm}[1]{\bigl\Vert#1\bigr\Vert} 
\newcommand{\Norm}[1]{\Bigl\Vert#1\Bigr\Vert} 
\newcommand{\Bnorm}[1]{\Bigl\Vert#1\Bigr\Vert} 
\newcommand{\Set}[2]{\{\,#1\,:\,#2\,\}}
\newcommand{\bset}[2]{\big\{\,#1\,;\,#2\,\big\}}
\newcommand{\Bset}[2]{\Big\{\,#1\,;\,#2\,\Big\}}
\newcommand{\field}[1]{\mathbb{#1}}
\newcommand{\N}{\field{N}}


\definecolor{dred}{rgb}{.8,0,0}
\definecolor{ddmagenta}{rgb}{0.7,0,0.9}
\definecolor{ddcyan}{rgb}{0,0.2,1.0}
\definecolor{Green}{rgb}{0.,0.5,0.}

\newcommand{\DDDS}{\color{red} }
\newcommand{\DDDE}{\color{black}}
\newcommand{\DelS}{\color{yellow} }
\newcommand{\DelE}{\color{black}}

%% file: cover.tex


\thispagestyle{empty} 

\phantom{1}

\vspace{0.5cm}
\noindent {\Large
{\bf Energy release and Griffith's criterion for phase-field fracture}
}

\vspace{24pt}

\vspace{24pt}

\begin{small}
{\bf E.~Maggiorelli and M.~Negri}

\medskip
{Department of Mathematics -  University of Pavia} 

{Via A.~Ferrata 1 - 27100 Pavia - Italy}

\vspace{36pt}
\noindent {\bf Abstract.} 
Phase field evolutions are obtained by means of time discrete schemes, providing (or selecting) at
each time step an equilibrium configuration of the system, which is usually computed by descent methods for
the free energy (e.g.~staggered and monolithic schemes) under a suitable irreversibility constraint on the
phase-field parameter. We study in detail the time continuous limits of these
evolutions considering monotonicity as irreversibility constraint and providing a general result, which holds
independently of the scheme employed in the incremental problem. In particular, we show that in the steady
state regime the limit evolution is simultaneous (in displacement and phase field parameter) and satisfies
Griffith’s criterion in terms of toughness and phase field energy release rate. In the unsteady regime the
limit evolution may instead depend on the adopted scheme and Griffith’s criterion may not hold. We prove also
the thermodynamical consistency of the monotonicity constraint over the whole evolution, and we study the
system of PDEs (actually, a weak variational inequality) in the steady state regime. Technically, the proof employs a suitable reparametrization of the time discrete points, whose Kuratowski limit characterizes the set of steady state propagation. The study of the quasi-static time continuous limit relies on the strong convergence of the phase field function together with the convergence of the power identity.

\bigskip
\noindent {\bf AMS Subject Classification: 49S05, 74A45.} 

\end{small}

%% file: notation.tex

\section{Notation \label{not}}

Assume that $\Omega$ is a bounded Lipschitz connected open set in 
$\mathbb{R}^2$. The boundary $\partial \Omega$ of $\Omega$ is split into $\partial_D \Omega$ (for Dirichlet boundary condition) and $\partial_N \Omega = \partial \Omega \setminus \partial \Omega_D$ (for Neumann boundary condition);  we assume that $\partial_D \Omega$ is relatively open in $\partial \Omega$. The spaces of the phase-field and of the displacement variables will be respectively
$$  \V := H^1 ( \Omega)\cap L^{\infty}(\Omega) ,  \qquad  \U := H^1 (\Omega , \R^2)  . $$
Although the phase-field functions $v \in \V$ may a priori take any real values, in the evolution we will have $\dot{v} \le 0$ (to model irreversibility) and $0 \leq v \leq 1$, with $v=1$ corresponding to no damage and  $v=0$ corresponding to maximum damage.  

The total energy takes the form 
$$
	\F (u, v) := \E  (u,v) + G_c \, \mathcal{L} (v) .
$$
In this work, we are not studying any limit as the internal length $\eps$ vanishes, so for convenience we omit the dependence on $\eps$ in the notation.

In the literature there are several options for the elastic energy $\E$ and the length term $\mathcal{L} $ that form $\F$, for a comprehensive view we refer to \cite{AmbatiGerasimovDeLorenzi_CM15,BourdFrancMar08,Wick20,VicentiniZolesiCarraraMauriniDeLorenzis_IJF24}. 
 Here, the elastic energy takes into account the volumetric-deviatoric split \cite{AmorMarigoMaurini_JMPS09,ChambolleContiFrancfort_ARMA18} and reads 
$$
	\E  (u ,v) := \int_\Omega W ( \strain (u) , v ) \, dx ,
$$
$$ 
	W ( \strain , v ) = \psi (v) W_+ ( \strain ) + W_- ( \strain )  , 
	\quad 
	W_+ (\strain) :=   \mu   | \strain_d |^2 + \kappa   | \strain_v^+ |^2 ,
	\quad
	W_- (\strain)  :=   \kappa   | \strain_v^- |^2,
$$%
where $\strain_v^+ : = \tfrac12 \mathrm{tr}^+ ( \strain) \boldsymbol{I}$ and $\strain_v^- : = \tfrac12 \mathrm{tr}^- ( \strain) \boldsymbol{I}$  are respectively the tensile and compressive components of the strain, and $\strain_{d} := \strain - \strain_v$ is the deviatoric component.  In \S \ref{GR} we will also employ a slightly different elastic energy to gain better properties of the energy release. The degradation function takes the form $\psi (v) :=  v^2 + \eta$ for $\eta >0$.

 With this energy split \cite{AmorMarigoMaurini_JMPS09} only the strain energy related to shear and expansion competes with the surface energy. 
Note that, the energy density $W (\strain ,v ):=\psi (v)W_+(\strain )+W_-(\strain) $ is differentiable with respect to $\strain $ and $\partial_{\strain} W ( \strain, v )$ gives the phase-field stress 
\begin{equation} \label{e.stress} \stress ( u,v) :=  2 \psi (v) \bigl( \mu \strain_d (u) + \kappa \strain^+_v (u) \bigr) - 2 \kappa \strain_v^- (u).  \end{equation}
 As far as the length term, we employ the classical quadratic functional
$$
	\mathcal{L} (v) := \tfrac12 \int_\Omega (v-1)^2+ | \nabla v |^2 dx,
$$
known as $AT_2$ in the engineering literature.

\medskip
Let $g (s) := \alpha (s) \hat g$ where $\hat g \in W^{1,q} (\Omega; \mathbb{R}^2) $ for $q>2$, and $\alpha\in C^1([0,1])$\bl, and define
$$  \U (s) : = \{  u \in H^1 (\Omega ; \R^2)  :   u = g (s) \mbox{ on  } \partial_D \Omega \}  \quad \text{for $s \in [0,1]$.} $$
Working in the rate-independent setting, $g$ provides the {\it loading path} for the system (independently of time); in the sequel we will find a suitable time parametrization, or control, of the loading path $g$ in such a way that the crack speed will remain finite. 
Due to the irreversibility constraint, the natural set of admissible variations 
for the phase field variable is   the cone $ \Xi = \{ \xi \in H^1 (\Omega) \cap L^{\infty } (\Omega): \xi \le 0 \}$, 
while the displacement variations are in the space  $  \Phi = \{  \phi \in H^1 (\Omega ; \mathbb{R}^2) : \phi = 0 \mbox{ on  } \partial_D \Omega \}$; thus 
equilibrium is characterized by 
\begin{align*}
	\partial_v \F ( u , v) [ \xi ] \ge 0 , \qquad \partial_u  \F ( u , v) [ \phi ] = 0,   \qquad 
	\text{for every $\xi \in \Xi$ and $\phi \in \Phi$.}
\end{align*}

\section{Statement of the main results}
\label{mainres}
In the study of rate-independent processes, it is customary to employ schemes based on a time discretization or on a discretization of the loading path. 
Given $v_0 \in \V$ with $0  \le v_0 \le 1$, let $u_0 \in \argmin \{  \F  (  u , v_0) : u \in \U (0) \} $. For $n\in \N$, $n>0$, set $s^n_k := k / n $ for $k=0,...,n$ and $v^n_0 : = v_0$. 
\bl Let us consider $( u^n_k, v^n_k)$ for $k=1,...,n$ such that
\begin{equation} \label{e.newsequence}
\begin{cases}
	u^n_{k} \in \argmin \{  \F  (  u , v^n_{k} ) : u \in \U (s_k^n)
\} 
\\
         v^n_{k} \in \argmin \{  \F ( u^n_{k} , v ) : v \in \V 
\mbox{ with } v \le v^n_{k-1} \} . 
\end{cases} 
\end{equation} 
In other terms, $(u^n_k, v^n_k)$ is a separate minimizer of the energy $\F$ under the constraints $u \in \U (s^n_k)$ and $v \le v^n_{k-1}$. 

\begin{example} \normalfont In practice, there are several ``sub-schemes'' to provide the update $(u^n_k , v^n_k)$: for instance, we may employ an energy descent method, such as staggered (see \S \ref{M-ele}), or global energy minimization. In the latter case,
$$
	(u^n_k, v^n_k ) \in \argmin \{  \F (u,v) : u \in \U (s^n_k) , \, v \in \V ,\, v \le v^n_{k-1} \} . 
$$
In the former, we introduce the auxiliary sequences $\{ u^n_{k,i} \}$ and $\{v^n_{k,i} \}$ defined by induction in the following way: $u^n_{k,0} = u^n_{k-1}$, $v^n_{k,0}  = v^n_{k-1}$, and then 
$$
	\begin{cases}
	u^n_{k,i+1} \in \argmin \{  \F  (  u , v^n_{k,i} ) : u \in \U (s_k^n)
\} 
\\
         v^n_{k,i+1} \in \argmin \{  \F ( u^n_{k,i+1} , v ) : v \in \V 
\mbox{ with } v \le v^n_{k-1} \} . 
\end{cases} 
$$
Then, see Proposition \ref{Elimt},  $v^n_{k,i} \to v_\infty$ in $L^2$ and $u^n_{k,i} \to u_\infty$ in $H^1$ (up to non-relabelled subsequences)  where $(u_\infty, v_\infty)$ satisfies \eqref{e.newsequence}; thus we set $u^n_k = u_\infty$ and $v^n_k = v_\infty$. 
Other choices could be constrained minimization \cite{RoerentropBoddinKneesMosler_CMAME24}, second order schemes \cite{LeonBaldelliMaurini_JMPS21} etc.

\end{example}

From \eqref{e.newsequence} by separate minimality we have 
\begin{align} \label{eqDISC}
	\partial_u  \F (u^n_{ k} ,v^n_{k} ) [ \phi ] = 0, \quad  \partial_v\F(u^n_{ k} ,v^n_{k} )[v-v^n_{k}]\geq 0,   \quad 
	\text{for every $\phi \in \Phi$ and $v\in \V$ with $v \le v^n_{k-1} $.}
\end{align}
As a straightforward consequence 
\begin{align} \label{eq0DISC}
\partial_v\F(u^n_{ k} ,v^n_{k} )[v^n_{k}-v^n_{k-1}  ] = 0.
\end{align}

In general, the ``speed''  $\| v^n_k - v^n_{k-1} \|_{H^1} / ( s^n_k - s^n_{k-1} )$ can be arbitrarily large; however, with a suitable time-parametrization of the loading path we can control it, as follows. 
\begin{lemma} \label{l.repar}  Let $\tilde{v}^n$ be the piecewise affine interpolation of $v^n_k$ in the points $s^n_k$. There exist a fixed time  $T$ and a sequence of uniformly Lipschitz, increasing  parametrizations $ c^n : [0,T] \to [0,1]$ such that the sequence $\{ v^n = \tilde{v}^n \circ c^n \}_{n \in \mathbb{N}}$ is bounded in $W^{1,\infty} ( 0,T ;  H^1 (\Omega))$.
Given $v^n$, let $u^n (t):=\argmin\{\E( u ,v^n(t)) : u\in \U( c^n(t))\}$; there exists $p>2$ such that  the sequence $\{ u^n \}_{n \in \mathbb{N}} $ is bounded in $W^{1,\infty} ( 0,T ; W^{1,p} (\Omega ; \mathbb{R}^2))$.
\end{lemma}

\noindent From the above result it follows that there exists a subsequence (non relabelled) such that 

\begin{itemize}\setlength\itemsep{0pt}
\item $c^n \weakstarto c$ in $W^{1,\infty} (0,T)$, and thus uniformly in $[0,T]$, 
\item $v^n \weakstarto v$ in $W^{1,\infty} ( 0, T ; H^1(\Omega) )$,
\item $u^n \weakstarto u$ in $W^{1,\infty} ( 0, T ; W^{1,p}(\Omega;\mathbb{R}^2))$ for some $p>2$. \bl
\end{itemize}
The parametrization $c : [0,T] \to [0,1]$ is Lipschitz continuous and increasing. Given the parametrization $c$ we will write by abuse of notation $\U(t) = \U ( c(t) )$ etc.  Clearly, $v$ and $u$ are Lipschitz continuous with the same constant of $v^n$ and $u^n$ respectively.

\begin{remark}\normalfont From the \textit{mechanical} point of view, the parametrizations $c^n$ and $c$ provide some sort of control of the boundary condition, ensuring finite speed of the crack, and actually arbitrarily small speed, which is consistent with the intuitive idea that quasi-static evolutions are slow in time; however, even if the evolution is reparametrized and the speed in controlled, in general we cannot avoid unstable branches, as explained in the sequel. 
From the \textit{mathematical} point of view, parametrizations allow to distinguish clearly steady and unsteady evolutions (see below) and to work with Lipschitz maps in Sobolev spaces (such a regularity would not be possible using only monotonicity and relying on Helly's Theorem for compactness).
\end{remark}

In order to study the limit evolution $(u,v)$ a couple of delicate properties are needed: strong convergence of $v^n$ and $u^n$ and a characterization of stable and unstable regimes (in the parametrized setting). The following result provides strong convergence, putting together Lemma \ref{ls.strong} and Corollary \ref{c.strong}.

\begin{lemma} \label{l.strong} 
Let $t \in [0,T]$, then $v^n (t) $ converge strongly to $v(t)$ in $H^1 (\Omega)$.
\bl Moreover, $u^n(t)$ converges to $u (t):=\argmin\{\E( u ,v (t))\,\,u\in \mathcal U(t)\}$ strongly in $W^{1,p} (\Omega ; \mathbb{R}^2)$ for some $p>2$.

\end{lemma}

In order to separate stable and unstable regimes we argue in the following way, inspired by the theory of visco-energetic evolutions \cite{MinottiSavare_ARMA18}. For every $n \in \mathbb{N}$, $n>0$, we define a suitable set of (parametrized) time discrete points $ \Upsilon^n:=\{t^n_{k} : k=0,...,n \} \subset [0,T]$, corresponding to the points $s^n_k$. There exists a subsequence (non relabelled) such that $\Upsilon^n$ converge to a set $I_s \subseteq [0,T]$ in the sense of Kuratowski. As we will see, the set $I_s$ is the set where the (parametrized) evolution is steady state (or stable), while the set $I_u = [0,T] \setminus I_s$ is the set where the evolution is unstable.  We will show that $I_s$ is a compact set of positive measure, 
see Corollary \ref{mis_pos}. 

The set $I_s$ depends on the time discrete sequence and in practice on the scheme used to provide it.  For example, global minimization algorithms may predict different evolutions from those obtained by energy descent methods, see e.g.~\cite{NegriOrtner_M3AS08,Maggiorelli24}: since the energy is non-convex it may happen that, in the incremental procedure, global minimization finds a configuration $(u^n_k, v^n_k)$ even if the potential energy is not decreasing along any path starting from $(u^n_{k-1}, v^n_{k-1})$. In this scenario, an energy descent algorithm would not predict $(u^n_k, v^n_k)$. Even within energy descent algorithms, alternate minimization and monolithic schemes may predict different evolutions, see e.g.~\cite{LeonBaldelliMaurini_JMPS21}.

\begin{remark} \normalfont Linear interpolation between $(u^n_{k-1}, v^n_{k-1})$ and $(u^n_k, v^n_k)$ is a very natural choice for our result, which is independent of the incremental sub-scheme. When the evolution is steady-state any interpolation would be fine, since $(u^n_{k-1}, v^n_{k-1})$ and $(u^n_{k}, v^n_{k})$ will converge to the same limit. However, in the unstable regime $(u^n_{k-1}, v^n_{k-1})$ and $(u^n_{k}, v^n_{k})$ converge to different configurations; in this case a precise description of the transition between them would require ad-hoc results depending on the specific incremental sub-scheme. 
\end{remark}

Next Theorem describes the behaviour of the evolution in the steady state regime by a sort of {\it evolutionary variational inequality} in Karush-Kuhn-Tucker form. 

\begin{theorem} \label{EVI}
 For every  $t \in I_s $ the limit configuration $(u(t), v(t))$ is an equilibrium point for the energy $\F (\cdot, \cdot)$, i.e., 
 \begin{align*}
	\partial_v \F ( u (t) , v (t) ) [ \xi ] \ge 0 , \qquad \partial_u  \F ( u (t) , v(t) ) [ \phi ] = 0,   \qquad 
	\text{for every $\xi \in \Xi$ and $\phi \in \Phi$.}
\end{align*}
Moreover, 
for a.e.~$t \in I_s$ it holds  
$$
	\partial_v \F ( u(t) ,  v (t) ) [ \dot{v}(t)] = 0
$$

\end{theorem} 

The statement of Theorem \ref{EVI} can be recast in terms of PDEs, under sufficient regularity, as stated in the next result; this is indeed the system of evolution equations employed in several engineering papers, see e.g.~\cite{VicentiniZolesiCarraraMauriniDeLorenzis_IJF24}.

\begin{theorem} \label{t.PDE} For $t \in I_s$ the solution $u$ solves the following PDE 
\begin{equation} \label{eq.PDEu}
	\begin{cases}
		\mathrm{div} ( \stress ( t ) ) = 0    & \Omega \\
		\stress ( t ) \, \hat{n} = 0 		& \partial_N \Omega \\
		u (t) = g (t) & \partial_D \Omega .
	\end{cases}
\end{equation}
where for simplicity $\stress(t)$ denotes the phase-field stress \eqref{e.stress}. 
Moreover, if $\dot{v} \in C (\bar\Omega)$ we have 
\begin{equation} \label{eq.PDEuu}
	\begin{cases}
		- G_c \Delta v (t) + G_c ( v (t) -1 ) + 2 v (t)  W_+ ( \strain(t) )   \le 0   & \Omega \\[2pt]
		\displaystyle \frac{\partial v (t) }{ \partial n} \le 0 	& \partial \Omega \\[2pt]
		\big[ - G_c \Delta v (t)+ G_c (v (t) -1) + 2 v (t)  W_+ ( \strain( t) )  \big] \dot{v} (t) = 0   & \Omega \\[2pt]
		\displaystyle \frac{\partial v (t) }{ \partial n} \dot{v} (t) = 0	& \partial \Omega 
	\end{cases}
\end{equation}
where $\strain (t) = \strain ( u(t))$. 
\end{theorem}

The limit evolution satisfies also the following energy identity. 

\begin{proposition} \label{p.pow-id}
For $a.e.~t \in [0,T]$ it holds
$$
	\F ( u(t) , v(t)) = \F ( u_0, v_0) + \int_{0}^t \mathcal{P} (s, u(s) , v(s) ) \, ds + \int_{I_u \cap (0,t)} \partial_v \F ( u(s) , v(s) ) [ \dot{v} (s) ] \, ds  , 
$$
where $\mathcal{P}$ is the power 
 \begin{align*} 
   \mathcal{P} (t, u(t),v(t))\, : = 
    \int_\Omega\stress(u(t),v(t)) : \strain(\dot{ a }(t) \hat g ) \, dx , 
\end{align*} 
and $a(t):=\alpha(c(t))$.  
\end{proposition}

\begin{remark}\label{cor.lip} \normalfont Observe that the function $t\mapsto (u(t),v(t)) \mapsto \F ( u(t) , v(t) )$ is Lipschitz continuous, as it is the composition of a Lipschitz continuous, bounded function with a locally Lipschitz function, (see Lemma \ref{lip}).  
If the evolution is steady state the set $I_u$ is empty and thus the classic energy identity holds. On the contrary, for unstable evolutions it may happen that 
$$   \int_{I_u \cap (0,t)} \partial_v \F ( u(s) , v(s) ) [ \dot{v} (s) ] \, ds \neq 0 . $$ 
\end{remark}

Finally, we want to show that the limit evolution satisfies Griffith's criterion. First, let us provide the notion of {\it (maximal) energy release} in the phase field context \cite{KneesNegri_M3AS17} 
\begin{align}
 \label{def.ERR}
	  \mathcal{G}_{\text{\sl max}}(t,v) &: = \mathrm{sup} \big\{ \! - \partial_v \tilde{\E} (t, v) 
[\xi] : \xi \in \Xi \text{ with } d \mathcal{L} (v) [\xi]  = 1  \big\} ,
\end{align} where $\tilde{\E} (t, v) $ is the reduced energy
$$
	\tilde{\E} ( t ,v) := \E (  u_{t,v} , v)
	\quad \mbox{ for } \quad
	u_{t,v} \in \argmin \{ \E ( u , v) : u \in \U(t)  \} .
$$
Other representations and properties of $\mathcal{G}_{\text{\sl max}}$ are contained in \S \ref{GR}. 
The energy release rate, as we defined it, is not a simple quantity to compute, since the derivative of the energy should be calculated in any admissible direction. Nevertheless, we will see in the sequel that there is no need to do so.  
The reason for introducing $\mathcal{G}_{\text{\sl max}}$ is to prove that the evolution $t \mapsto v(t)$ satisfies Griffith criterion, with this notion of energy release; this result is summarized in next Theorem, where, in order to facilitate its statement, we denote $\ell(t) := \mathcal{L} ( v(t) )$. 

\begin{theorem} \label{t.KKT} \bl
Assume that $v_0 \not \equiv 1$. 
The limit evolution $v$ satisfies $0 \le v (t)  \le 1$ for every $t \in [0,T]$, the irreversibility constraint $\dot{v} (t) \le 0$ and the thermodynamic consistency condition $\dot{\ell} ( t ) =d\mathcal{L}(v(t))[\dot{v}(t)] \ge 0$ a.e.~in $[0,T]$. The following Griffith's criterion with maximal energy release rate holds. In the steady state regime we have 
\begin{itemize}
\item $  \mathcal{G}_{\text{\sl max}}( t  , v(t) ) \le G_c$ everywhere~in $I_s$, 
\item $(  \mathcal{G}_{\text{\sl max}}( t , v(t) ) - G_c ) \, \dot{\ell} (t) = 0$ a.e.~in $I_s$,
\item $\partial_v \tilde\E ( t , v )  [ \dot{v}(t)]  = -  \dot{\ell}(t) \, \mathcal{G}_{\text{\sl max}} ( t ,v(t) )$ a.e.~in $I_s$. 
\end{itemize}
In the unsteady regime we have instead
\begin{itemize}
\item $  \mathcal{G}_{\text{\sl max}}(t, v(t)) \ge G_c$ on a subset of positive measure of $I_u$.
\end{itemize}
\end{theorem}

\begin{remark} \normalfont The irreversibility of the crack is modelled by the monotonicity of the phase field variable $v$. This hypothesis does not directly imply the monotonicity of the dissipated energy, which is however true for the evolution, as stated in the above Theorem \ref{t.KKT}. Note that in the steady state regime, the above properties correspond precisely to \eqref{e.GR1}-\eqref{e.GR2} in the sharp crack setting.  \end{remark}

%% file: preliminaries.tex

\section{Energy variations and continuous dependence \label{prel}}

This brief section is devoted to some preliminary results that will be widely used in the following, starting from some properties of the energy density $W$. 
Let $\mathbb{M}^2_{\text{\sl sym}}$ be the space of symmetric (real) $2 \times 2$ matrices.
For some of the proosf of next Lemmata, we refer to \cite{AlmiNegri_ARMA20,HM}.

\begin{lemma}\label{lemmaAN}The energy density $W:\mathbb{M}^2_{\text{\sl sym}} \times\mathbb{R}\to \mathbb{R}_+$ is of class $C^{1,1}_{\text{\sl loc}}$. Moreover, there exist two  constants $c, \,C>0$ such that, for every $v\in[0,1]$ and every $\strain_1,\strain_2\in\mathbb{M}^2_{\text{\sl sym}}$, the following holds:
\begin{enumerate}
\item $\bigl(\partial_{\strain} W(\strain_1,v)-\partial_{\strain} W(\strain_2,v)\bigr):\bigl(\strain_1-\strain_2\bigr)\geq c|\strain_1-\strain_2|^2$,
\item $ |\partial_{\strain} W(\strain_1,v)-\partial_{\strain} W(\strain_2,v) |\leq C|\strain_1-\strain_2|$. 
\end{enumerate}
Observe that, since $\partial_{\strain} W(0,v)=0$, for every $\strain\in\mathbb{M}^2_{\text{\sl sym}}$ we have $|\partial_{\strain} W(\strain,v) |\leq C|\strain|.$
\end{lemma}

\begin{lemma} \label{l.deriv} The functional $\F (  \cdot, \cdot)$ is sequentially weakly lower semicontinuous in $\U (s) \times \V$. 
Moreover,  for $v \in \V $ 
and $u \in \U (s)$ the partial derivative of $\F$ w.r.t.~$u$ is given by
\begin{gather}
   \partial_u \F (u,v) [\phi ] = \partial_u \E (u,v) [\phi ]  =  \int_\Omega \partial_{\strain} W ( \strain (u) , v) : \strain (\phi) \, dx = \int_\Omega  \stress(u,v):\strain(\phi) \, dx ,
\end{gather}
for every $\phi \in \Phi$. If $v  \in \V$ and $u \in \U$, then for every $\xi \in \Xi $,
\begin{gather}
  \partial_v \F (u,v) \, [\xi]  = 
   \int_\Omega  2v\xi \, W_+ (\strain(u)) \, dx + G_c \int_\Omega (v-1) \xi + \nabla v 
\cdot \nabla \xi  \, dx.
\label{e.dzF}
\end{gather}

\end{lemma}

 \begin{corollary}\label{esp_quad}
 There exists $C>0$ such that, for every $s \in [0,1]$ and every $v_0,\,v_1\in\V$, it holds:
$$|\E(u_1,v_0)-\E(u_0,v_0)|\leq C\|u_1-u_0\|_{H^1}^2,$$
 where $u_0=\argmin\{\E(u,v_0):\,u\in\U (s) \}$ and $u_1 \in \U (s)$. 
 \end{corollary}
 \proof
 Let us consider the convex combinations $u_r=ru_1+(1-r)u_0$ where $r \in [0,1]$ and write:
 $$\E(u_1,v_0)-\E(u_0,v_0)=\int_0^1 \frac{d}{dr}\E(u_r,v_0)\, dr=\int_0^1 \partial_u\E(u_r,v_0)[u'_r]\, dr,$$ where $u'_r = u_1 - u_0$. 
 Since $\partial_u\E(u_0,v_0)[\phi]=0$ for every $\phi \in \Phi$, we can continue with:
 \begin{align*}\int_0^1 \partial_u\E(u_r,v_0)[u'_r]\, dr&=\int_0^1 \partial_u\E(u_r,v_0)[u'_r]\, dr-\partial_u\E(u_0,v_0)[u'_r]\, dr\\&=\int_0^1 \langle\stress(u_r,v_0)-\stress(u_0,v_0),\strain(u'_r)\rangle_{L^2}\, dr.
 \end{align*}
 It follows that
 $$|\E(u_1,v_0)-\E(u_0,v_0)|\leq \int_0^1 \|\stress(u_r,v_0)-\stress(u_0,v_0)\|_{L^2}\|\strain(u'_r)\|_{L^2}\, dr.$$
 By Lemma \ref{lemmaAN}, we get:
 $$\|\stress(u_r,v_0)-\stress(u_0,v_0)\|_{L^2}\leq C\|\strain(u_r)-\strain(u_0)\|_{L^2}= Cr\|\strain(u_1-u_0)\|_{L^2}\leq C \|u_1-u_0\|_{H^1}.$$
Since $\|\strain(u'_r)\|_{L^2}=\|\strain(u_1-u_0)\|_{L^2}\leq C\|u_1-u_0\|_{H^1}$, the proof is concluded.
 \qed

\medskip
A higher integrability of the strains is often needed and it is obtained from the fact that $g(s) \in W^{1,q} ( \Omega; \mathbb{R}^2)$ for $q>2$.

\begin{lemma} \label{l.reg} There exists $\tilde p \in (2,q]$ and $C>0$ such that for every $v \in \V$ with $\| v\|_{L^\infty } \le 1$ and every $s \in [0,1]$ it holds
$$
	\| u \|_{W^{1, \tilde{p}}} \le C , 
$$
where $u = \argmin\{\E(u,v): \,u \in\U (s) \}$. 
\end{lemma}

In the next two lemmata we provide two results of continuous dependence. 
 
\begin{lemma} \label{le.KRZ2.5} 
Let $\tilde p$ be as in Lemma \ref{l.reg}. 
For $i=1,2$, let $s_i \in [0,1]$ and $v_i \in \V$ with $0\le v_i \le 1$ and denote by $u_i$ the minimizer of $\E ( \cdot , v_{i})$ over \,$\U(s_i)$.  
Then 
$$
\| u_2 - u_1 \|_{W^{1, p}}  \le C \bigl(| 
s_2 - s_1 | + \| v_2 - v_1 \|_{L^r}  \bigr)
$$
where $C>0$ is independent of $s_i$ and $v_i$ while 
$\frac{1}{r}=\frac{1}{p}-\frac{1}{\tilde{p}}$ for $p\in[2,\tilde{p}]$.
\end{lemma}  

\begin{remark} \label{r.pr} \normalfont
 Let us fix for the rest of the paper $\tilde p \in ( 2 , q]$ 
as the exponent obtained in  Lemma \ref{l.reg}. Moreover, let us also fix $p\in(2,\tilde{p})$ and accordingly $r \in [1, +\infty)$ 
as in Lemma \ref{le.KRZ2.5}. 
\end{remark}

\begin{lemma} \label{lip} The energy $\F$ is locally Lipschitz continuous in $W^{1,p}(\Omega,\mathbb{R}^2)\times \V $. 
\end{lemma}
\proof
We consider the elastic and the dissipated energy separately. Regarding the latter, $ \mathcal{L} (v)=\tfrac12 \norm{v-1}_{H^1}^2$  
is locally Lipschitz continuous since the $H^1$-norm is Lipschitz. As for the elastic energy, for every $(u_1,v_1),\,(u_2,v_2)\in W^{1,p}(\Omega,\mathbb{R}^2)\times H^1(\Omega)$, it holds: 
\begin{align*}
  |\E(u_2,v_2)-\E(u_1,v_1)|
    &\leq\int_{\Omega}\big|W(\strain(u_2),v_1)-W(\strain(u_1),v_1)\big|+\big|W(\strain(u_2),v_2)-W(\strain(u_2),v_1) \big| \, dx \\
    &=\int_{\Omega}\big|(v_1^2+\eta)\bigl(W_+(\strain(u_2))-W_+(\strain(u_1))\bigr)+  \big(W_-( \strain (u_2) ) - W_- (\strain (u_1)) \big)\big|\\ &\quad+\int_{\Omega} | v_2^2-v_1^2 | W_+(\strain(u_2)) \, dx
    \end{align*}
Using the inequality $| a^\pm - b^\pm | \le | a - b|$ we easily get
$$
	 \big| (a^\pm)^2 - (b^\pm)^2 \big| = | a^\pm - b^\pm | ( a^\pm + b^\pm) \le | a - b | ( a^\pm + b^\pm) .
$$
Hence we can continue with 
    \begin{align*}
   |\E(u_2,v_2) & -\E(u_1,v_1)|  \le \\ &\leq\int_{\Omega}(v_1^2+\eta)\bigl(\mu|\strain_d(u_2-u_1)|\big(|\strain_d(u_2)|+|\strain_d(u_1)|\big) +\kappa|\strain_v(u_2-u_1)|\big(\strain_v^+(u_2)+\strain_v^+(u_1)\big)\bigr) \, dx 
   \\ &\quad+\int_\Omega \kappa|\strain_v(u_2-u_1)|\big(\strain_v^-(u_2)+\strain_v^-(u_1)\big)\\ &\quad+\int_{\Omega}\big|(v_2-v_1)(v_2+v_1)\big(\mu|\strain_d(u_2)|^2+\kappa|\strain_v^+(u_2)|^2\big)\big|dx
   \\ &\leq C  \|v_1^2+\eta\|_{L^\infty} \|u_2-u_1\|_{W^{1,p}}  \big(\|u_2\|_{W^{1,p'}}+\|u_1\|_{W^{1,p'}}\big) 
   \\[3pt] & \quad +C \|v_2-v_1\|_{L^{(p/2)'}}  \|v_2+v_1\|_{L^\infty}   \|u_2\|^2_{W^{1,p/2}}   .
\end{align*} 
From the continuous embedding of $H^1 (\Omega)$ in $L^r (\Omega)$ for every $1\leq r<+\infty$, 
it follows that
\begin{align*}
    |\E(u_2,v_2)-\E(u_1,v_1)|\leq C\bigl(\norm{u_2-u_1}_{W^{1,p}}+\norm{v_2-v_1}_{H^1}\bigr)
\end{align*} 
and the lemma is proven. \qed

Next Lemma follows from \cite{KneesRossiZanini_M3AS13}, together with Lemma \ref{l.reg}. 

\begin{lemma} \label{l.Frechet} The reduced elastic energy $\tilde{\E} ( t, \cdot)$ is Fr\'{e}chet differentiable in $\V$, with respect to the $H^1$-norm. 
\end{lemma}

\begin{lemma}
\label{lem_dk_est_z_1WEAK} 
Let 
 $u_1,u_2\in W^{1,p}(\Omega;\mathbb{R}^2)$, $ v_0\in \calV$ with $0\leq  v_0\leq 1$ and
\begin{align}
v_1&=\argmin\Set{\calF(u_1,v)}{ v \in \calV, v\leq v_0}, \label{min_constr1} \\
v_2&=\argmin\Set{\calF(u_2,v)}{v\in \calV, v\leq v_0}.  \label{min_constr2} 
\end{align}
Then 
\begin{align}
\label{est_dk_1}
\hspace{-6pt}\norm{v_1- v_2}_{H^1}
\leq C \norm{u_1- u_2}_{H^1} ,
\end{align}
where the constant $C$ is independent of $u_i$ and $v_i$.
\end{lemma}

\proof
We adapt the argument of \cite[Proposition 1]{AlmiNegri_ARMA20}. The solutions of the constrained minimization problems \eqref{min_constr1} and \eqref{min_constr2} verify the following inequality:
$$\partial_v \F (u_i,v_i) \, [v-v_i] \geq0 \quad \text{
for all $v\in\V$ such that $v\leq v_0$, for $i=1,2$}.$$ 
Therefore, for $i=2$ and taking $v=v_1\leq v_0$, we obtain that $\partial_v \F (u_2,v_2) \, [v_2-v_1]=-\partial_v \F (u_2,v_2) \, [v_1-v_2] \leq0 $. Moreover,  for $i=1$ and setting $v=v_2\leq v_0$, $\partial_v \F (u_1,v_1) \, [v_2-v_1] \geq0 $. Hence
\begin{align*}
\langle \partial_v\calF(u_2, v_2), v_2 - v_1\rangle \leq  \langle \partial_v\calF(u_1, v_1), v_2 - v_1\rangle ,
\end{align*}
which yields  
\begin{align}
\langle \partial_v\calF(u_2, v_2) -  \partial_v\calF(u_2, v_1), v_2 - v_1\rangle
 \le \langle \partial_v\calF(u_1, v_1) -  \partial_v\calF(u_2, v_1), v_2 - v_1\rangle.
\label{est_dk_4}
\end{align}
The left hand side of \eqref{est_dk_4} is greater than $G_c \norm{v_2-v_1}_{H^1}^2$. Indeed,
\begin{align*}
\langle \partial_v\calF(u_2, v_2) & -  \partial_v\calF(u_2, v_1), v_2 - v_1\rangle \\ & = 
   \int_\Omega v_2  (v_2-v_1)\, W_+ (\strain(u_2))  \bl \, dx +G_c \int_\Omega (v_2-1) (v_2-v_1)+ \nabla v_2 
\cdot \nabla (v_2-v_1) \, dx \,\\
&\quad-\biggl(\int_\Omega v_1 (v_2-v_1)\, W_+ (\strain(u_2))   \, dx +G_c \int_\Omega (v_1-1) (v_2-v_1)+ \nabla v_1 
\cdot \nabla (v_2-v_1) \, dx \biggr)\\
&= \int_\Omega (v_2-v_1)^2 W_+ (\strain(u_2))  \, dx +G_c \int_\Omega (v_2-v_1)^2+|\nabla (v_2-v_1) |^2\, dx \,\geq G_c  \norm{v_2-v_1}_{H^1}^2.
\end{align*}
As for the right hand side of \eqref{est_dk_4}, we obtain:
\begin{align*}
\langle \partial_v\calF(u_1, v_1)  & -  \partial_v\calF(u_2, v_1), v_2 - v_1\rangle \\ & = 
\int_\Omega v_1 (v_2-v_1)\, W_+ (\strain(u_1)) \, dx +G_c \int_\Omega (v_1-1) (v_2-v_1)+ \nabla v_1
\cdot \nabla (v_2-v_1) \, dx \,\\
&\quad-\biggl(\int_\Omega v_1 (v_2-v_1)\, W_+ (\strain(u_2))  \, dx +G_c \int_\Omega (v_1-1) (v_2-v_1)+ \nabla v_1 
\cdot \nabla (v_2-v_1) \, dx \biggr) \\
&=\int_\Omega v_1(v_2-v_1) \big( W_+ (\strain(u_1)) - W_+ (\strain(u_2))  \big)\, dx \\ 
& \leq C \int_\Omega v_1 | v_2-v_1 | | \strain(u_1-u_2) | \big(|\strain(u_1)|+|\strain(u_2)|\big)dx\\
& \leq C \norm{v_1}_{L^{\infty}} \norm{v_2-v_1}_{L^r} 
\norm{u_1-u_2}_{H^1}\big(\|u_1\|_{W^{1,p}}+\|u_2\|_{W^{1,p}}\big)\\ & \leq C \norm{v_2-v_1}_{H^1}  \norm{u_1-u_2}_{ H^1} 
, 
\end{align*}
where we 
applied H\"older inequality 
with $\tfrac1r + \tfrac12 + \tfrac1p =1 $ 
and exploited the continuous embedding of $H^1(\Omega)$ in $L^r(\Omega)$ for any $1 \le r < \infty$. This finishes the proof.
\qed

%% file: conv-time.tex

\section{Convergence of time discrete evolutions}

\subsection{Control in time and interpolation} 

We present a preliminary result that will then be used in the proof of Lemma \ref{l.repar}: 
\begin{lemma} \label{l.length} There exists a constant $\tilde C$ such that $\sum_{k=1}^{n}\|v_k^n-v_{k-1}^n\|_{H^1}\leq \tilde C$ \bl for every $n \in \mathbb{N}$. 
\end{lemma}

\bl
\proof
\bl By separate minimality, $v^n_{k-1} \in \argmin \{  \F ( u^n_{k-1} , v ) : v \in \V  \mbox{ with } v \le v^n_{k-2} \}$ and therefore, we also have $v^n_{k-1} \in \argmin \{  \F ( u^n_{k-1} , v ) : v \in \V  \mbox{ with } v \le v^n_{k-1} \}$, since $v^n_{k-1} \le v^n_{k-2}$. 
By definition $v^n_{k} \in \argmin \{  \F ( u^n_{k} , v ) : v \in \V \mbox{ with } v \le v^n_{k-1} \}$, hence,  Lemma \ref{lem_dk_est_z_1WEAK} yields $\|v^n_k-v^n_{k-1}\|_{H^1}\leq  C \| u^n_k-u^n_{k-1}\|_{ H^1}$. Moreover, by minimality of $u^n_k$, Lemma \ref{le.KRZ2.5} provides $\| u^n_k-u^n_{k-1} \|_{W^{1,p}}\leq C (|s^n_k-s^n_{k-1}|+\|v^n_k-v^n_{k-1}\|_{L^r})$. 
Putting together these estimates, we get that for $k=1,\dots,n$ \bl
\begin{align}\label{vk}\|v^n_k-v^n_{k-1}\|_{H^1}\leq  C \|u^n_k-u^n_{k-1}\|_{W^{1,p}}\leq C (|s^n_k-s^n_{k-1}|+\|v^n_k-v^n_{k-1}\|_{L^r}) .
\end{align}

Let $\theta\in (0,1)$. Then, for $s:=\tfrac{r(1-\theta)}{1-r\theta}\in(r,+\infty)$ 
and $\delta>0$, to be chosen a posteriori, by interpolation in $L^p$-spaces we get: \bl
\begin{align*}\|v^n_k-v^n_{k-1}\|_{L^r}&\leq  \delta ^{1 - \theta} \|v^n_k-v^n_{k-1}\|_{L^s}^{1-\theta} \ \delta^{\theta-1}\|v^n_k-v^n_{k-1}\|_{L^1}^{\theta} \\
&\leq (1-\theta)\delta \|v^n_k-v^n_{k-1}\|_{L^s}+  \theta \delta^{1-1/\theta}\|v^n_k-v^n_{k-1}\|_{L^1} , 
\end{align*}
where in last inequality we used  Young's inequality. 
Therefore it follows from \eqref{vk} and from the embedding of $H^1(\Omega)$ in  $L^s(\Omega)$ that
\begin{align*}\|v^n_k-v^n_{k-1}\|_{H^1}\leq C(\tfrac1n+ C_1 \delta \|v^n_k-v^n_{k-1}\|_{H^1}+ C_2 \delta^{1 - 1/\theta} \|v^n_k-v^n_{k-1}\|_{L^1})
\end{align*} and if we take $\delta$ sufficiently small (in such a way that $c :=  C C_1 \delta <1$) we can absorb the second term in the right-hand side into the left-hand side.
 Taking the sum over $k$, by the monotonicity of $v^n_k$ with respect to $k$ we finally get  
\begin{align*}(1-c)\sum_{k=1}^n\|v^n_k-v^n_{k-1}\|_{H^1}\leq \sum_{k=1}^nC (\tfrac1n+C'_2 \|v^n_k-v^n_{k-1}\|_{L^1}) \leq  C + C' \|v^n_n-v^n_{0}\|_{L^1} , 
\end{align*} hence the thesis follows.
\qed

{\bf Proof of Lemma \ref{l.repar}.} 
 Let $S>0$. First of all we define the auxiliary non-uniform steps $\rho^n_k$ according to the length $\|v^n_{k}-v^n_{k-1}\|_{H^1}$ at the given step $s^n_k$, i.e., for $k=1,\dots,n$:
\begin{equation} \label{e.taunk}
\rho^n_k= \mathrm{max} \{   \tfrac1n , \tfrac1S \|  v^n_{k}-v^n_{k-1}\|_{H^1}  \}  = 
\begin{cases}
\tfrac1n & \text{if } \|v^n_{k}-v^n_{k-1}\|_{H^1} \leq\tfrac{S}n \\
\tfrac1S \|v^n_{k}-v^n_{k-1}\|_{H^1}  & \text{otherwise. } 
\end{cases}
\end{equation}
Note that $\rho^n_k \ge 1/n$. We define $r^n_0:=0$, and for $k=1,\dots,n$,\, $r_k^n:=r^n_{k-1} + \rho^n_k = \sum_{j=1}^k\rho^n_j$. It holds that $T^n:=r^n_n$ is finite: indeed, by Lemma \ref{l.length} 
$$T^n=\sum_{k=1}^n\rho^n_k\leq \sum_{k=1}^n\bigl(\tfrac1n+ \tfrac1S \|v^n_{k}-v^n_{k-1}\|_{H^1} \bigr)\leq 1+\tfrac{\tilde C}{S}.$$
 Note that $T^n \ge 1$ since $\rho^n_k\geq 1/n$. 
Let $w^n : [0,T^n] \to [0,1]$ be the piecewise affine map which interpolates $s^n_k$ in $r^n_k$. Clearly $w^n$ is monotone  increasing  and bounded. Moreover, being $\rho^n_k \ge 1/n$, we have that for $r\in [r^n_{k-1}, r^n_k ]$:
\begin{equation} \label{1lip}
 |\dot{w}^n(r)|=\frac{s^n_{k} - s^n_{k-1}}{r^n_k - r^n_{k-1}} = \frac{s^n_{k} - s^n_{k-1}}{\rho^n_k} \le n  (s^n_{k} - s^n_{k-1}) = 1 ,    
\end{equation}
and thus $w^n$ is $1$-Lipschitz. 
There exists a \bl (non-relabelled) \bl subsequence of $T^n$ converging to a certain $T \ge 1$. We reparametrize each of the \bl corresponding (non-relabelled) parametrizations \bl $w^{n}$ on the interval $[0, T]$ as follows: $c^{n}(t):=w^{n}(t \tfrac{T^{n}}{T})$.  In this way  $c^{n}(0)=0$ and $c^{n}(T)=1$. Moreover,  the parametrizations $c^n$ result to be uniformly Lipschitz, since:
\begin{equation} \label{c.lip}
    \big|\dot{c}^{n}(t)\big|=\big|\dot{w}^{n}(t \tfrac{T^{n}}{T})\tfrac{T^{n}}{T}\big|\leq \sup \{T^{n}\}=:C'.
\end{equation}
Set $t^{n}_k:=r^{n}_k\tfrac{T}{T^{n}}$ and $\tau_k^{n}:= t^{n}_k-t^{n}_{k-1}=\rho^{n}_k\tfrac{T}{T^{n}}$, we define the piecewise affine interpolate of $v^{n}_k$ in the times $t^{n}_k$:
\begin{equation}
    v^{n}(t):=v_{k-1}^{n}+(t-t_{k-1}^{n})\frac{v_k^{n}-v_{k-1}^{n}}{\tau_k^{n}}    \quad\text{ for $t\in[t_{k-1}^{n},t_k^{n})$  } 
\end{equation}
and observe that it has controlled speed. Indeed, for $t\in[t_{k-1}^{n},t_k^{n})$
\begin{equation*}
    \|\dot{v}^{n}(t)\|_{H^1}=\frac{\norm{v^{n}_k-v_{k-1}^{n}}_{H^1}}{\tau^{n}_k},
\end{equation*} 
and if $\|v^{n}_{k}-v^{n}_{k-1}\|_{H^1} \leq \tfrac{S}{{n}}$, then $\tau^{n}_k=\tfrac{1}{{n}}\tfrac{T}{T^{n}}$, and thus 
$  \|\dot{v}^{n}(t)\|_{H^1}\leq \tfrac{S/{n}}{1/{n}}\tfrac{T^{n}}{T}\leq SC'$. Otherwise, $\|v^{n}_{k}-v^{n}_{k-1}\|_{H^1} =S\rho^{n}_k$ and $$\|\dot{v}^{n}(t)\|_{H^1}= \displaystyle \frac{\norm{v^{n}_k-v_{k-1}^{n}}_{H^1}}{\tau^{n}_k}=\frac{S\rho^{n}_k}{\rho^{n}_k T /T^{n}}=S \tfrac{T^{n}}{T}\leq SC' . $$ 
This finishes the first part of the proof. 
 Let us now check that the sequence $\{u^n\}_{n \in \N }$ is uniformly bounded in $W^{1,\infty}( 0,T ; W^{1, p}(\Omega ;\mathbb{R}^2))$. For $t_1, t_2 \in [0,T]$, by minimality and Lemma \ref{le.KRZ2.5} 
\begin{align*}         
	\norm{ u^n(t_2) - u^n(t_1) }_{W^{1, p }}
	& \le
	C  | c^n (t_2)  - c^n (t_1) | + C  \norm{v^n (t_2) - v^n (t_1)}_{L^{r}}.
\end{align*}
This, by Sobolev embedding of $H^1(\Omega)$ in $L^r(\Omega)$,  
 concludes the proof. \qed

\bl
\begin{remark} \normalfont
Clearly there exists a subsequence (not relabelled) of $c^n$ that converges uniformly in $[0,T]$. 
    The limit parametrization $c : [0,T] \to [0,1]$ is 
Lipschitz continuous, increasing, and surjective.
 As anticipated in Section \ref{mainres}, given the parametrization $c$ we will write, by abuse of notation, $\U(t) = \U ( c(t) )$ etc.

Differently from \cite{KneesNegri_M3AS17} we do not employ a reparametrization with respect to the arc length (in the intrinsic norms), rather we provide a time parametrization of the loading path. Note that, thanks to the control $c^n$ introduced in Lemma \ref{l.repar}, the evolutions are Lipschitz continuous in time. 
\end{remark}

\subsection{Compactness }

\begin{lemma} \label{l.weakLIP} The piecewise affine interpolates $v^n$ 
 converge weakly* (up to non-relabelled subsequences) to a limit $v$ in $W^{1,\infty}( 0,T ;H^1(\Omega))$; therefore $v^n$ converge to $v$ in $L^r(\Omega)$  for all $r<+\infty$ uniformly in $[0,T]$. Moreover, $u^n$ converge weakly* to a limit $u$ in 
 $ W^{1,\infty}( 0,T ; W^{1, p}(\Omega ; \mathbb{R}^2))$. 
\end{lemma}

\proof
By Lemma \ref{l.repar}, the piecewise affine interpolate $v^n$ weakly* converges  (up to non-relabelled subsequences) to an evolution $v$ in $W^{1,\infty}([0,T];H^1(\Omega))$, that, by Aubin-Lions Lemma, is compactly embedded in $C([0,T];L^r(\Omega))$ for all $r<+\infty$. Therefore $v^n $ converges to $v$ uniformly in $L^r(\Omega)$.
In the same way,  $u^n$ weakly* converges (up to non-relabelled subsequences) to a limit $u$ in $ W^{1,\infty}([0,T];W^{1,p}(\Omega, \mathbb{R}^2))$.  
 \qed

Next, we will prove strong convergence in $H^1(\Omega)$ of the phase-field variable. Note that this property does not follow from Aubin-Lions lemma. Moreover, it cannot follow solely from unilateral minimality, since at the moment we know that the obstacles $v^n_k$ converge weakly. The proof follows instead from unilateral minimality together with a (delicate) construction of  back-in-time admissible competitors, inspired by \cite[Proposition 2.13]{KimuraNegri_NoDEA21}. 
To this purpose, we must first introduce a partition of the time set $[0,T]$ as follows.
For every $n \in \mathbb{N}$, we define $ \Upsilon^n:=\{t^n_{k} : k=0,...,n \} \subset [0,T]$. 
We call $I_s\subseteq [0,T]$ the Kuratowski limit of this sequence of sets (up to non-relabelled subsequences). This is a compact set and therefore $I_u:=[0,T]\setminus I_s=\bigcup_{j\in \mathbb{N} }(a_j,b_j)$ is the countable union of disjoint open intervals.  For a brief overview on Kuratowski convergence, we refer to Appendix \ref{Kur}. 

The following Lemma immediately follows from Lemma \ref{le.KRZ2.5} and will be widely employed in the sequel. 

\begin{lemma} \label{ele} Let $p$ and $r$ be as in Remark \ref{r.pr}. If $t^m \to t$ and $v^m$ converges to $v$ in $L^r(\Omega)$, then $u^m:=\argmin\{\F(u\,,v^m): u \in \U(t^m )\}$ converges to $u \in \argmin \{  \F  (u \, , v ) : u \in \U(t)\}$ in \,$W^{1, p }(\Omega\ \mathbb{R}^2)$. 
\end{lemma}
We are now ready to present the main result on the strong convergence of the sequences $\{v^n(t)\}_{n\in \mathbb{N}}$ and $\{u^n(t)\}_{n\in \mathbb{N}}$ as they are defined in Lemma \ref{l.repar}. 

\begin{lemma} \label{ls.strong} Let $p$ and $r$ be as in Remark \ref{r.pr}. 
Let $t \in I_s$, then $v^n (t) $ converge strongly to $v(t)$ in $H^1 (\Omega)$ and $u^n(t)$ 
converges to $u(t)$ strongly in $W^{1, p} (\Omega ; \mathbb{R}^2)$. 
\end{lemma}

\proof For sake of simplicity we  divide the proof in three 
parts: we start by showing estimate \eqref{e.bj}, that will be used, in the second part, to prove that $v^n(t^n_{k_n})$ strongly converges to $v(t)$ in $H^1 (\Omega)$, for every $t\in I_s$ and every $\{t^n_{k_n}\}_{n\in\mathbb{N}}\in \Upsilon^n$ converging to $t$.   Finally, in the third part of the proof, we will prove the strong convergence of $v^n(t)$ and $u^n(t)$. 

{\bf I.} 
For $1< s<\frac{\tilde p}{2}$ let $\mathcal{W} = \{ W_+^{s} ( \strain (u^n_k) ) : n \in \mathbb{N}, \, 0 \le k \le n \}$. By Lemma \ref{l.reg}  and by de la Vall\'ee-Poussin Lemma we know that the family $\mathcal{W}$ is equi-integrable: 
indeed, being $u \in W^{1,\tilde p}(\Omega) $, it holds that for $s<\frac{\tilde p}{2}$, there exists $c > 1$ such that  $\int_\Omega W_+^{c s }( \strain (u^n_{k}) ) \,  dx\leq C$. 

Therefore, there exists a modulus of continuity $\psi : \mathbb{R}^+ \to \mathbb{R}^+$  such that $\lim_{x \to 0^+} \psi(x) = 0$ and $\int_A W_+^{s } ( \strain (u^n_{k}) ) \, dx \le \psi ( | A |)$ for every measurable set $A \subset \Omega$ and every choice of $n \in \mathbb{N}$ and $0 \le k \le n $. 

We claim that there exists $C>0$ such that the following holds: let $z \in \V$ be such that $z \le v_{0}$, then for every $n \in \mathbb{N}$ and every index $0 \le  m \le n$ we have
\begin{equation} \label{e.bj}
      \int_{\{ v^n_m < z \}} | \nabla v^n_m  |^2 \, dx \le  \int_{\{ v^n_m < z \}} | \nabla z |^2 \, dx + C \psi^{1/s}  \big( | \{ v^n_m < z \} |   \big)  . 
\end{equation}

 Let $j := \max \{ k : v^n_k  \ge z \}$ (remember that $ v^n_0 = v_0 \ge z $ and thus $j \ge 0$).  If $m \le j$ then the set $\{ v^n_m < z \}$ is empty and there is nothing to prove in \eqref{e.bj}. 
Otherwise, we provide a construction of admissible competitors going back-in-time up to $t^n_j$. 
For $j \le k \le m$ let $A^n_k = \{ v^n_k < z \}$. Note that $A^n_j = \emptyset$, moreover, the monotonicity of $v^n_k$ (with respect to $k$) implies that the sets $A^n_k$ are increasing (with respect to $k$); in particular the family $\{ A^n_k \setminus A^n_{k-1} : \text{ for } k = j  +1, ..., m \}$ is a disjoint partition of the set $A^n_m = \{  v^n_m < z \}$, which appears in \eqref{e.bj}. 
\begin{itemize}
\item
Let $k= j +1$. We introduce the auxiliary function $z^n_k  = \max \{ z , v^n_k \}$. Note that $z^n_k \le v^n_{k-1}$ since  $v^n_k \le v^n_{k-1}$ (by monotonicity) and  $z  \le v^n_{k-1}  = v^n_j$ (by definition of $j$). Hence $z^n_k$ is an admissible competitor for the incremental problem and by minimality we get 
$\F ( u^n_k , v^n_k ) \le \F ( u^n_k ,  z_k^n  )$, i.e., 
\begin{align*}
	\int_\Omega ((v^n_k)^2 + \eta ) W_+\bl ( \strain(u^n_k) ) \, dx  & + \tfrac12 G_c \int_\Omega (v^n_k -1)^2 + | \nabla v^n_k |^2 \, dx \\ & \le  
	\int_\Omega  ((z^n_k)^2 + \eta ) W_+\bl ( \strain(u^n_k) )  \, dx +  \tfrac12 G_c \int_\Omega (z^n_k-1)^2 + | \nabla z^n_k |^2 \, dx .
\end{align*}
As $ v^n_k \le z^n_k $ we get $(v^n_k-1)^2 \ge (z^n_k-1)^2$ and thus 
\begin{align*}
	\int_\Omega ((v^n_k)^2 + \eta ) W_+\bl ( \strain(u^n_k) ) \, dx  & + \tfrac12 G_c \int_\Omega | \nabla v^n_k |^2 \, dx \\\
	& \le \int_\Omega((z^n_k)^2 + \eta ) W_+\bl ( \strain(u^n_k) ) \, dx  + \tfrac12 G_c \int_\Omega | \nabla z^n_k |^2 \, dx . 
\end{align*}

 Note that $\{ z^n_k \neq v^n_k \} =  \{ v^n_k < z \} = A^n_k \subset A^n_m$. Therefore, we can ``localize'' the above inequality in $A^n_k$, where $z^n_k=z$, obtaining 
\begin{align} \label{e.11}
     \tfrac12 G_c \int_{A^n_k} | \nabla v^n_k |^2 \, dx 
      &\le 
 \tfrac12 G_c \int_{A^n_k} | \nabla z |^2 \, dx  
     +  \int_{A^n_k}   \big( z^2  -  ( v^n_k )^2  \big) W_+ ( \strain ( u^n_k )) \, dx  \nonumber \\
     & \le 
     \tfrac12 G_c \int_{A^n_k} | \nabla z |^2 \, dx  
     +  \int_{A^n_m}   \big( (v^n_{k-1})^2  -  ( v^n_k )^2  \big) W_+ ( \strain ( u^n_k ))  \, dx \nonumber \\
          & \le 
      \tfrac12 G_c \int_{A^n_k} | \nabla z |^2 \, dx  
     + \tfrac12 C \| v^n_{k-1}  -  v^n_k \|_{L^{s'}} \|  W_+ ( \strain ( u^n_k )) \|_{L^{ s} (A^n_m)}  \nonumber \\ 
     & \le \tfrac12 G_c 
      \int_{A^n_k} | \nabla z |^2 \, dx  
     + C \| v^n_{k-1}  -  v^n_k \|_{L^{s'}} \psi^{1/s} \bl ( | A^n_m | ) . 
\end{align} 

\item For $j+1 < k \le m$ we introduce the auxiliary functions  
$$z^n_k = \min \{ \max \{ z , v^n_k \} , v^n_{k-1} \} . $$
Note that 
\begin{equation}   \label{e.wnkins}
	\begin{cases}  
		v^n_k \le v^n_{k-1} < z  & \text{in $A^n_{k-1}$} , \\
		v^n_k < z \le  v^n_{k-1}   & \text{in $A^n_k \setminus A^n_{k-1}$} , \\
		z \le v^n_k \le v^n_{k-1}   & \text{in $\Omega \setminus A^n_k$} .	
	\end{cases}
\end{equation} 
Hence the function $z^n_k$ can be equivalently written as 
\begin{equation}  \label{e.wnk}
	z^n_k = 
	\begin{cases}  
		v^n_{k-1} & \text{in $A^n_{k-1}$,} \\
		z & \text{in $A^n_k \setminus A^n_{k-1}$,} \\
		v^n_k & \text{in $\Omega \setminus A^n_k$.} 
	\end{cases}
\end{equation}
By construction $v^n_k \le z^n_k \le v^n_{k-1}$, hence $z^n_k$ is again an admissible competitor in the incremental problem and by minimality we get first $ \F ( u^n_k , v^n_k ) \le \F ( u^n_k , z_k^n  ) $ and then (arguing as above) 
\begin{align*}
	\int_\Omega  ((v^n_k)^2 + \eta) W_+\bl ( \strain(u^n_k) ) \, dx  & + \tfrac12 G_c \int_\Omega | \nabla v^n_k |^2 \, dx 
	\\&\quad\le  
	\int_\Omega ( (z^n_k)^2 +\eta) W_+\bl ( \strain(u^n_k) ) \, dx  + \tfrac12 G_c \int_\Omega | \nabla z^n_k |^2 \, dx  .
\end{align*}
Using \eqref{e.wnk}, the ``localization'' of the previous inequality in $A^n_k\subset A_m^n$ becomes
\begin{align*}
     \tfrac12 G_c \int_{A^n_k}  | \nabla v^n_k |^2 \, dx 
      & \le   
      \tfrac12 G_c \int_{A^n_{k-1}} | \nabla v^n_{k-1} |^2 \, dx  +  \tfrac12 G_c \int_{A^n_k \setminus A^n_{k-1}} | \nabla z |^2 \, dx  
     \\&\quad+  \int_{A^n_k}  \big( (v^n_{k-1})^2 -  (v^n_k)^2 \big) W_+\bl ( \strain ( u^n_k )) \, dx   
      \\
      & \le \tfrac12 G_c \int_{A^n_{k-1}} | \nabla v^n_{k-1} |^2 \, dx  +  \tfrac12 G_c \int_{A^n_k \setminus A^n_{k-1}} | \nabla z |^2 \, dx 
      \\&\quad
      + C \| v^n_{k-1}  -  v^n_k \|_{L^{s'}} \psi^{ 1/s} ( | A^n_m | ) . 
\end{align*}
\end{itemize}
Note that the above estimate actually holds also for $k=j+1$, indeed $A^n_j = \emptyset$ and thus it is equivalent to \eqref{e.11}.  Hence for every $j +1 \le k \le m$ we have 
\begin{align*}
	\int_{A^n_k}  | \nabla v^n_k |^2 \, dx -  \int_{A^n_{k-1}} | \nabla v^n_{k-1} |^2 \, dx  & \le   \int_{A^n_k \setminus A^n_{k-1}} | \nabla z |^2 \, dx + C \| v^n_{k-1}  -  v^n_k \|_{L^{s'}} \psi^{1/s} ( | A^n_m | )  . 
\end{align*}
Let us take the sum  for $j+1 \le k \le m$. Remember that: $\{ A^n_k \setminus A^n_{k-1} \}$ provide a disjoint partition of $A^n_m = \{ v^n_m < z \}$, $ A_j^{n} = \emptyset$, and that $v^n$ is of bounded variation in $L^r (\Omega)$ for any $1 \le r < +\infty$, since it is uniformly Lipschitz continuous in $H^1(\Omega)$. Hence we get \eqref{e.bj}. 

{\bf II.} 
Since $t \in I_s$, there exists a sequence $\{t^n_{k_n}\}_{n\in \mathbb{N}} $ such that $t^n_{k_n} \to t$.  Even though each index $k_n$ obviously depends on $n$, we chose to omit this dependence for simplicity of notation, setting $k:=k_n$. 
As $v^n \weakstarto v$ in $W^{1,\infty} (0,T ; H^1(\Omega))$ we have $v^n_{k} = v^n (t^n_{k}) \weakto v(t)$ in $H^1(\Omega)$; by Sobolev embedding $v^n_{k} \to v (t)$ in $L^r(\Omega)$ for every $r \in [1,+\infty)$ and thus $v^n_{k} \to v(t)$ in measure, i.e., for every  $\delta >0 $ 
\begin{equation} \label{mes}
	| \{ | v^n_{k} - v(t) | > \delta  \} | \to 0. 
\end{equation}

By minimality we have $u^n_k \in \argmin \{ \E ( u , v^n_k ) : u \in \U (t^n_k) \}$.  
Since $v^n_k \to v(t)$ in $L^r(\Omega)$ and $t^n_k \to t$, by Lemma \ref{ele}, $u^n_k \to u (t)$ strongly in $W^{1,p} ( \Omega ; \mathbb{R}^2)$ where  $ p=\tfrac{r\tilde p}{\tilde p+r}$ and
\begin{equation}
\label{forte}
	u (t) \in \argmin \{ \E ( u , v (t) ) : u \in \U (t) \} . 
\end{equation}
Let us show that $v^n_{k} \to v(t)$ strongly in $H^1(\Omega)$. By lower semicontinuity of the norm we have 
$$
	\int_\Omega | \nabla v (t) |^2 \, dx \le \liminf_{n \to +\infty} \int_\Omega | \nabla v^n_{k} |^2 \, dx , 
$$
therefore it is enough to show that 
$$
	\int_\Omega | \nabla v (t) |^2 \, dx \ge \limsup_{n \to +\infty} \int_\Omega | \nabla v^n_{k} |^2 \, dx ,
$$
 so that $ \lim_{n \to +\infty} \int_\Omega | \nabla v^n_{k} |^2 \, dx=\int_\Omega | \nabla v (t) |^2 \, dx $. Indeed, combining the convergence in 
norm with the fact that $v^n_{k}\weakto v(t)$ in $H^1(\Omega)$, we get strong convergence in $H^1(\Omega)$. 
\medskip
\\For $\delta >0$ let $z := \mathrm{max} \{ 0 , v(t) - \delta\} $ and $z^- := \min \{ z , v^n_k \}$. Let us also introduce the set 
$$   Z^n_k = \{ z^- < z \} = \{ v^n_k < z =  \mathrm{max} \{ 0 , v(t) -\delta \}  \} =  \{ v^n_k < v(t) -\delta  \}	.	$$
Note that $| Z^n_k | \to 0$ by convergence in measure \eqref{mes}. 
Moreover, since $z^- \le v^n_k \le v^n_{k-1}$, it turns out that $z^-$ is a competitor for the problem 
$$
	\min \{ \F ( v , u^n_k ) : v \in \V, \,  v \le v^n_{k-1} \} 
$$
which is solved by $v^n_k$. Therefore we have  $ \F ( v^n_k , u^n_k) \le \F ( z^- , u^n_k)$, i.e.,
\begin{align*}
	\int_\Omega ( (v^n_k)^2 + \eta ) W_+\bl( \strain(u^n_k) ) & \, dx   + \tfrac12 G_c \int_\Omega (v^n_k -1)^2 + | \nabla v^n_k |^2 \, dx \\ & \le  
	\int_\Omega  ((z^-)^2 + \eta ) W_+\bl ( \strain(u^n_k) ) \, dx  + \tfrac12 G_c \int_\Omega (z^--1)^2 + | \nabla z^- |^2 \, dx .
\end{align*}
Since $v^n_k \ge z^-$ we can write 
\begin{align*}
	\int_\Omega (v^n_k -1)^2 + | \nabla v^n_k |^2 \, dx \le  \int_\Omega (z^--1)^2 + | \nabla z^- |^2 \, dx .
\end{align*}
 In the set $Z^n_k$ we have $v^n_k < z$ and thus $z^- = v^n_k$.
Hence, we can ``localize'' the previous integral inequality in the set $\Omega \setminus Z^n_k = \{ z^- \neq v^n_k \} = \{ z^- = z \}$ and re-arranging the terms we get 
\begin{align*}
	\int_{\Omega \setminus Z^n_k} | \nabla v^n_k |^2 \, dx 
	& \le 
	 \int_{\Omega \setminus Z^n_k} | \nabla z^-  |^2 \, dx + \int_{\Omega \setminus Z^n_k}   (z^- -1)^2 - (v^n_k -1)^2 \, dx \\
	 & =   
	 \int_{\Omega \setminus Z^n_k} | \nabla z  |^2 \, dx + \int_{\Omega \setminus Z^n_k}   ( z -1)^2 - (v^n_k -1)^2 \, dx.
\end{align*}
Note that $Z^n_k = \{ v^n_k < z \}$, hence by \eqref{e.bj} with $m=k$ we get 
$$
    \int_{Z^n_k} | \nabla v^n_k  |^2 \, dx \le  \int_{Z^n_k} | \nabla z  |^2 \, dx + C \psi^{1/s}  \big( | Z^n_k |  \big)   . 
$$
Taking the sum of the previous inequalities and remembering that $\nabla z = \nabla v(t)$ provides 
\begin{align*}
     \int_{\Omega} | \nabla v^n_k  |^2 \, dx \le  \int_{\Omega} | \nabla v(t) |^2 \, dx + \int_{\Omega \setminus Z^n_k}   ( z -1)^2 - (v^n_k -1)^2 \, dx  + C \psi^{1/s}  \big( | Z^n_k |   \big)  . 
\end{align*}
Clearly last term is infinitesimal. For the second term we write 
$$
	|  ( 1- z )^2 - (1 - v^n_k )^2 | = | ( v^n_k - z )   (  2 - z - v^n_k ) | \le C |  v^n_k - z  | \le C | v^n_k - v(t) | + C \delta . 
$$
Hence
$$
	\limsup_{n \to \infty} \int_{\Omega \setminus Z^n_k}   ( z -1)^2 - (v^n_k -1)^2 \, dx  
	\le 
	\limsup_{n \to \infty}  \int_\Omega | v^n_k - v(t) | \, dx + C \delta | \Omega | \le C' \delta . 
$$
We conclude by the arbitrariness of $\delta > 0$. 

{\bf III.} In conclusion, applying Lemma \ref{l.repar}, we get:
$$
	\| v^n (t) - v(t) \|_{H^1} \le \|  v^n (t) - v^n (t^n_{k}) \|_{H^1} + \|  v^n (t^n_{k}) - v(t) \|_{H^1}  \le C | t^n_{k} - t | + \| v^n_{k} - v(t) \|_{H^1} \to 0 . 
$$
The strong convergence of $u^n (t)$ follows by Lemma \ref{le.KRZ2.5}. \qed

 \begin{remark} \normalfont 
We recall that the set $I_u$, being the complement of a compact subset of $[0,T]$, is the (at most) countable union of its connected components, denoted $(a_j, b_j)$ for $j \in \mathbb{N}$. Therefore, the set of right or left-isolated points $I_{\text{\sl isol}}= \left( \bigcup_{j \in \mathbb{N}} a_j  \right) \cup \left( \bigcup_{j \in \mathbb{N}}  b_j \right)$ has zero measure. 
\end{remark}

 We will now show that on each component of $I_u$ the limit parametrization $c$ is constant and the limit evolution $v$ is affine.  To do so, let us state some preliminary results.
\begin{lemma}  Let's consider the connected component $(a ,b )$ of $I_u$. There exists a sequence $\{t^n_{k_n} \in \Upsilon^n\}_{n \in \mathbb{N}}$ such that $t_{k_n}^n \to a$ and $t_{k_n +1}^n\to b$.
\label{conv_ajbj}
\end{lemma}

\proof 
First of all, observe that for every $\eps>0$, there exists $N_\eps\in \mathbb{N}$ such that  $\Upsilon^n\cap (a +\eps,b -\eps)=\emptyset$ for all $n>N_\eps$. 
Indeed, if there was a sequence $\{ t^{n_i}_{k_i} \}_{i\in \mathbb{N}} $ such that $t^{n_i}_{k_i}\in \Upsilon^{n_i}\cap  (a +\eps,b -\eps)$, we could take a 
subsequence converging to a limit  $t \in  [a +\eps,b -\eps]$; this $t$ would be in Kuratowski limit superior of $\Upsilon^n$, that is exactly $I_s$. However, by definition $[a +\eps,b -\eps] \subset I_u = [0,T] \setminus I_s$. 

Now, set $a < a' < b$  and define $a_n :=\max \{t^n_{k} \in \Upsilon^{n}\,\,:\,\,t^n_{k}\leq a' \}$.  This sequence must converge to $a$; 
let's prove it by contradiction.
If there was a sequence  $\{ n_i \}_{i\in \mathbb{N}}$ and $\eps > 0$ (sufficiently small) such that  $a_{n_i} > a +\eps$ then 
$ a_{n_i}\in(  \Upsilon^{n_i} \cap ( a +\eps ,a'] )\subset ( \Upsilon^{n_i}\cap(a +\eps,b -\eps) )$ which is empty for all $n>N_\eps$.
On the other hand, if there was a sequence $\{ n_i \}_{i\in \mathbb{N}} $  such that $ a_{n_i} <a -\eps$, by definition of $a_n$, 
$\Upsilon^{n_i}\cap [a -\eps, a' )=\emptyset$; hence, there wouldn't exist elements $  t^{n_i}_{k_i} \in \Upsilon^{n_i}$ 
converging (as $i\to \infty$) to $a $,  that, however, belongs to the Kuratowski limit of $\Upsilon^n$.

As we have seen, $\Upsilon^n \cap (a +\eps,b -\eps)=\emptyset$ for $n\geq N_\eps$.  Let $k_n$ be such that $ t^{n}_{k_n}=a_n$, then 
$t_{k_n+1}^n\geq b -\eps$ at least definitely. In particular, for $n\geq N_\eps$, $t_{k_n+1}^n=\min\{t^n_{k}\in \Upsilon^{n}\,\,:\,\,t^n_{k}\geq b -\eps\}$. 
We thus proceed as we did above and obtain that $t^n_{k_n+1}\to b$.
\qed
\begin{lemma}\label{ccost}The function $c$ is constant on each of the connected components of $I_u$.
\end{lemma}
\proof
Let's consider the connected component $(a ,b )$ of $I_u$. From Lemma \ref{conv_ajbj} we know that  there exists a sequence $\{t^n_{k_n}\in \Upsilon^n\}_{n\in \mathbb{N}}$ such that $t_{k_n}^n\to a $ and $t_{k_n+1}^n\to b $. 
 Now, since $c^n\to c$ uniformly in $[0,T]$, it holds that $s^n_{k_n}=c^n(t^n_{k_n})\to c(a )$ and $s^n_{k_n+1}=c^n(t^n_{k_n+1})\to c(b )$. Hence $c(b) - c(a) = \lim_{n \to \infty} s^n_{k_n+1} - s^n_{k_n} = \lim_{n \to \infty} 1/n = 0$, 
and the Lemma is proven. 
\qed

\begin{corollary}\label{mis_pos}
    The set $I_s$ has positive measure.
\end{corollary}
\proof
By Lemma \ref{ccost}, 
$\dot{c} =0$ in the open set $I_u$ 
therefore  $\text{supp}(\dot c) \subset I_u^c = I_s$. Now, the set $\text{supp}(\dot c)$ has positive measure, indeed,
\begin{align*}
    \int_{(0,T)}\dot{c}(t)\,dt=c(T)-c(0)=1
\end{align*}
and, being $0\leq \dot c \leq 1$,
\begin{align*}
    1=\int_{(0,T)}\dot{c}(t)\,dt= \int_{\text{supp}(\dot c)}\dot{c} (t) \,dt\leq |\text{supp}(\dot c)|.
\end{align*}
The proof is concluded.  \qed

\begin{proposition} Let $(a,b)$ be a connected component of the set $I_u$. 
The limit evolution $v$ is the affine interpolate between $v(a)$ and $v(b)$.
\label{v.aff}
\end{proposition}
\proof
Let $t\in(a ,b )$ and write it as $t=\theta b +(1-\theta)a $ for a certain $\theta \in (0,1)$. Let $t^n_{k_n}$ be as in Lemma \ref{conv_ajbj} and, as above, let us write $t^n_k$ for simplicity. 
As shown in step II of the proof of Lemma \ref{ls.strong}, 
$v^n(t^n_{k})\to v(a )$ and $v^n(t^n_{k+1})\to v(b )$ strongly in $H^1(\Omega)$. For $n$ sufficiently large, $t\in(t^n_{k}, t^n_{k+1})$, so we can write it as $t=\theta^n t^n_{k+1}+(1-\theta^n)t^n_{k}$, where 
$$\theta^n=\frac{t-t^n_{k}}{t^n_{k+1}-t^n_{k}} \to  \frac{t-a }{b -a }=\theta.$$
By definition 
$$ v^n(t)=\theta^nv^n(t^n_{k+1})+(1-\theta^n)v^n(t^n_{k})\to \theta v(b )+(1-\theta)v(a )$$
strongly in $H^1(\Omega)$, which is indeed the affine interpolate between $v(a )$ and $v(b )$.
\qed

From Proposition \ref{v.aff} together with Lemma \ref{ls.strong} we get the following result. 

\begin{corollary} \label{c.strong} Let $t \in [0,T] \setminus I_s$, then $v^n(t)$ converge strongly to $v(t)$ in $H^1(\Omega)$ and 
$u^n(t)$ converge strongly to $u(t)$ in $W^{1, p} (\Omega ; \mathbb{R}^2)$. 
\end{corollary}
This Corollary, together with Lemma \ref{ls.strong} gives Lemma \ref{l.strong}. 

\subsection{Equilibrium}

\begin{theorem}
\label{EQ} For every  $t \in I_s$ the limit configuration $(u(t), v(t))$ is an equilibrium point for the energy $\F (\cdot, \cdot)$, which  means that
\begin{equation}\label{equ} 
     \partial_u \F (u(t),v(t)) [\phi ]=0
	\quad \text{for all $\phi\in \Phi$,} \qquad \, \,
    \partial_v \F (u(t),v(t)) \, [\xi] \geq0 \,
\quad \text{for all $\xi \in \Xi.$}
\end{equation}
\end{theorem}
\proof
By definition of Kuratowski limit, there exist a sequence $\{t^n_{k_n} \in \Upsilon^n\}_{n\in \mathbb{N}}$ such that  $t^n_{k_n} \to t$. 
As usual, we write for simplicity $t^n_k$. \bl 
We have seen in the proof of Lemma \ref{ls.strong} that $v^n (t^n_{k}) \to v(t)$ in $H^1(\Omega)$, and thus in $L^r(\Omega)$ for any $r < \infty$,  
and that $u^n (t^n_{k}) \to u(t)= \argmin \{  \E  ( u , v (t) ) : u \in \U(t)\}$ in $W^{1, p} (\Omega ;  \mathbb{R}^2 )$ for $p \in (2,\tilde{p})$;  
 in light of these results, from \eqref{eqDISC} we obtain that 
\begin{equation}\label{eq1}
     \partial_u \F (u(t),v(t)) [\phi ]=  \int_\Omega \stress(u(t),v(t))  : \strain(\phi) \, dx  = 0
	\qquad \text{for all $\phi\in \V.$} 
\end{equation}
By definition  $v^n(t^n_{k}) \in \argmin \{  \F ( 
u^n_{k} , v ) : v \in \V 
\mbox{ with } v \le v^n_{k-1} \}$, and hence:
\begin{equation}\label{ineq_disc}
\int_\Omega v^n(t^n_{k})\xi \, W_+(\strain(u^n_k) ) \, dx + G_c \int_\Omega (v^n(t^n_{k})-1) \xi + \nabla v^n(t^n_{k})
\cdot \nabla \xi  \, dx\geq0 \,
\qquad \text{for all $\xi \in \Xi.$}
\end{equation}
As 
aforementioned $v^n(t^n_{k})$ converges to $v(t)$ strongly in 
 $H^1(\Omega)$ and $ W_+ ( \strain(u^n_{k} )) \to W_+(\strain(u(t))$ in $L^{\frac{p}{2}} (\Omega)$ where  $\tfrac{p}{2} > 1$.
Thereby, 
passing to the limit in \eqref{ineq_disc}, we get: 
\begin{equation} \label{eq2}
 \partial_v \F (u(t),v(t)) \, [\xi]  =\int_\Omega v (t) \xi \, W_+ (\strain( u(t))   \, dx + G_c \int_\Omega (v(t)-1) \xi + \nabla v (t)
\cdot \nabla \xi  \, dx\geq0 \,
\end{equation}
for all $\xi \in \Xi$. \qed

\subsection{Power identity}

In this section we will prove Proposition \ref{p.pow-id}, i.e., for every $t \in [0,T]$ 
$$
	\F ( u(t) , v(t)) = \F ( u_0, v_0) + \int_{0}^t \mathcal{P}  (s, u(s) , v(s) ) \, ds + \int_{I_u \cap (0,t)} \partial_v \F ( u(s) , v(s) ) [ \dot{v} (s)) ] \, ds  .
$$
We recall the definition of power (of external forces)
 \begin{align} \label{pext}
    \mathcal{P} (t, u(t),v(t))\, : = 
    \int_\Omega \stress(u(t),v(t)) : \strain(\dot{ a }(t) \hat g ) \, dx,     
\end{align} 
where $a(t):=\alpha(c(t))$.  
Notice that this is a good definition for the power: indeed, being $u(t)$ a minimum for the energy, by Theorem \ref{t.PDE} and by Green formula,  if $u$ is sufficiently regular it holds 
\begin{align*}
    \mathcal{P} (t, u(t),v(t)) & = 
    \phantom{|}_{H^{-1/2} ( \partial_D \Omega)} \langle \stress(u(t),v(t))\cdot n , \dot{a} (t) \hat g  \rangle_{H^{1/2}(\partial_D\Omega)} \nonumber  \\
   & = \int_{\partial_D \Omega} (\stress(u(t),v(t))\cdot n) \cdot \gamma(\dot{a} (t) \hat g ) \, d s  ,
\end{align*} 
this is indeed the product between force and the velocity imposed on the boundary. 
Finally note that  
\begin{equation}
     \mathcal{P} (t, u(t),v(t))  = \partial_u \F ( u(t) , v(t) ) [ \dot{u}(t) ]  \label{e.pow-wow}, 
\end{equation}  
since by Lemma \ref{l.repar} we can write  $\dot u(t)$ as the sum of $\dot g(t) = \dot a(t) \hat g$ and a function $w(t) \in \Phi $. 
 \begin{theorem} \label{t.enbal} For a.e.~$t \in I_s, $ it holds 
  \begin{align*} 
  \dot{\F} ( u (t) , v(t) ) =  \mathcal{P} (t, u(t), v(t)) . 
\end{align*}  
\end{theorem}

\proof
By Theorem \ref{EQ}, $(u(t), v(t))$ is an equilibrium point for the total energy $\F (\cdot, \cdot)$, and since $\dot{v}(t)$ is in $H^1(\Omega)$ and $\dot{v}(t)\leq 0$, by \eqref{equ} the following estimate for the total derivative of $\F$ holds almost everywhere in $I_s$: 
\begin{equation}\label{ineq}
\dot{ \F} (u(t),v(t))=\partial_u \F ( u(t) ,  v (t) ) [ \dot{u} (t)]+\partial_v \F ( u(t) ,  v(t) ) [ \dot{v}(t)]\geq \mathcal{P} (t, u(t),v(t)) ,
\end{equation}  
where we used \eqref{e.pow-wow}. We recall Remark \ref{cor.lip}, where we proved that $t \mapsto \F ( u(t) , v(t))$ is Lipschitz continuous, and thus a.e.~differentiable. To derive the thesis it is now sufficient to prove the inverse inequality.

{\bf I.} 
In this first part of the proof we ``go back'' to the time discrete setting, in order to prove that 
\begin{align}    \label{e.kanbis}
     \F (u^n_{k+1} , v^n_{k+1}) & \le   \F(u^n_{k}, v^n_{k})+\int_{t_{k}^n}^{t_{k+1}^n}\P^{n} (t,u^n (t) , v^{n} (t) ) \,dt + C | t^n_{k+1} - t^n_k |^2 ,
\end{align}
where  
$$	\P^{n} (t,u, v) := \int_\Omega \stress(u,v) :  \strain( \dot a^n(t) \hat g )  \, dx . $$ 
 We set 
$a^n_{k}:= \alpha ( s^n_k )= \alpha \circ c^n ( t^n_k) $  for $k=0,...,n$ and let $a^n$ represent the piecewise affine interpolant of $a^n_{k}$ in $t^n_k$. 
By minimality $\F (u^n_{k+1} , v^n_{k+1}) \le  \F(u^n_{k+1}, v^n_{k}) = \E (u^n_{k+1} , v^n_{k}) + G_c  \mathcal{L} (  v^n_{k})$. We define for any $j=0,...,n,$ $\hat u_j^n:=\argmin\{\E(u,v_j^n)\,:\,u=\hat g\,\text{ on }\partial_D\Omega\}$ and note that, by linearity, we can express $\E (u^n_{k+1} , v^n_{k})$ as $(a^n_{k+1})^2\E (\hat u^n_{k+1} , v^n_{k}) $.  According to Corollary \ref{esp_quad}, $\E (\hat u^n_{k+1} , v^n_{k})\le \E (\hat u^n_{k} , v^n_{k})+C\|\hat u_{k+1}^n-\hat u^n_k\|^2_{H^1}$ and by Lemma \ref{le.KRZ2.5}, last term is bounded by $C\|v_{k+1}^n-v^n_k\|^2_{H^1}\leq C|t^n_{k+1}-t^n_k|^2$, with the inequality following from Lemma \ref{l.repar}. Summarizing, 
\begin{align}\label{in1}\E ( u^n_{k+1} , v^n_{k})\le (a^n_{k+1})^2\E (\hat u^n_{k} , v^n_{k})+C'|t^n_{k+1}-t^n_k|^2.
\end{align}
We write the first addendum as
\begin{align*}(a^n_{k+1})^2\E (\hat u^n_{k} , v^n_{k})&=(a^n_{k})^2\E (\hat u^n_{k} , v^n_{k})+(a^n_{k+1}-a^n_{k}) \left( 2a^n_k+(a^n_{k+1}-a^n_{k})\right) \E (\hat u^n_{k} , v^n_{k}) \\
&=\E ( u^n_{k} , v^n_{k})+\bigg(2a^n_k(a^n_{k+1}-a^n_{k})+(a^n_{k+1}-a^n_{k})^2\bigg)\E (\hat u^n_{k} , v^n_{k})\\
&\leq \E ( u^n_{k} , v^n_{k})+2a^n_k(a^n_{k+1}-a^n_{k})\E (\hat u^n_{k} , v^n_{k})+C|t^n_{k+1}-t^n_k|^2 ,
\end{align*}
where the last inequality comes from the fact that $\E (\hat u^n_{k} , v^n_{k})$ is bounded and $a^n$ is Lipschitz continuous. Observe that, by definition of $\hat u^n_k$, $\E ( \hat u^n_{k} , v^n_{k})=\tfrac12\int_\Omega\stress(\hat u^n_k,v^n_k):\strain(\hat g)dx$, since, being $\hat u^n_k-\hat g$ an admissible variation, $\partial_u\E ( u^n_{k} , v^n_{k})[\hat u^n_k-\hat g]=0$. By bilinearity we hence get that $$2a^n_k(a^n_{k+1}-a^n_{k})\E (\hat u^n_{k} , v^n_{k})=\int_\Omega\stress( u^n_k,v^n_k):\strain((a^n_{k+1}-a^n_{k+1}) \hat g) \, dx$$ and \eqref{in1}  can finally be rewritten as:
\begin{align}\label{Pnk}\E ( u^n_{k+1} , v^n_{k})\le \E (u^n_{k} , v^n_{k})+\int_{t^n_k}^{t^n_{k+1}}\P^n(t,u^n_k,v^n_k) \, dt+C|t^n_{k+1}-t^n_k|^2.
\end{align}
Now, in order to obtain \eqref{e.kanbis}, 
we add and subtract $\P^{n} (t,u^n(t), v^n(t))$ to $\P^n(t,u^n_k,v^n_k)$: 
\begin{align*}
\int_\Omega\bigg(\stress(u^n_k,v^n_k)-\stress(u^n(t),v^n(t))\bigg):\strain(\dot a^n(t)\hat g)dx+\P^{n} (t,u^n(t), v^n(t)).
\end{align*}
We estimate the first term as follows (writing $\pm$ to add and subtract):
\begin{align}\label{first}
\int_\Omega\bigg(\stress(u^n_k,v^n_k)-\stress(u^n(t),v^n(t))\bigg):\strain(\dot a^n(t)\hat g)dx\pm \int_\Omega\stress(u^n(t),v^n_k):\strain(\dot a^n(t)\hat g) \, dx
\\=\dot a^n(t)\bigg[\int_\Omega\bigg(\stress(u^n_k,v^n_k)-\stress(u^n(t),v^n_k)\bigg):\strain(\hat g) \, dx \nonumber
\\+ \int_\Omega\bigg(\stress(u^n(t),v^n_k)-\stress(u^n(t),v^n(t))\bigg):\strain(\hat g) \, dx\bigg].\nonumber
\end{align}
By estimate 2 of Lemma \ref{lemmaAN}, the first addendum is lower than $C\|u^n_k-u^n(t)\|_{H^1}$.
As for the second addendum, we rewrite it explicitly:
\begin{align*}\int_\Omega\big((v^n_k)^2-v^n(t)^2\big) & \big(\mu \strain_d(u^n(t))+\kappa \strain_v^+(u^n(t))\big):\strain(\hat g) \, dx
\\ & \leq\|v^n_k-v^n(t)\|_{L^{s}}\|v^n_k+v^n(t)\|_{L^\infty}\|u^n(t)\|_{W^{1,p}}\|\hat g\|_{W^{1,q}},
\end{align*}
where $\frac{1}{s}+\frac1p+\frac1q=1.$
By the continuous embedding of $H^1(\Omega)$ in $L^{s}(\Omega)$,
\begin{equation*}\eqref{first}\le C \big(\| u^n_k - u^n (t) \|_{  H^1}+\| v^n_k - v^n (t) \|_{  H^1 \bl }\big)  
\le C' | t - t^n_k |  ,
\end{equation*}
where we used the fact that $u^n$ and $v^n$ are Lipschitz continuous (see Lemma \ref{l.repar}).
As a consequence,
$$
     \int_{t_{k}^n}^{t_{k+1}^n}\P^n(t,u^n_k,v^n_k)
     \, dt \le \int_{t_{k}^n}^{t_{k+1}^n}\P^{n} (t,u^n (t) , v^{n} (t) ) \,dt +C' | t^n_{k+1} - t^n_k |^2 
$$
and \eqref{e.kanbis} is proven.

{\bf II.} 
Starting from the estimate obtained in the discrete setting, let us now see what happens in the limit. Since the set of isolated points of $I_s$ has zero measure, we can state that for a.e. $t\in I_s$,  there exists a sequence $\{t_j\in I_s\}_{j\in \mathbb{N}}$ such that $t_j \to t^+$ or $t_j \to t^-$ with $t_j \neq t$. Let us consider the first case (we use the same reasoning for the second). Let $t^n_l, t_h^n\in \Upsilon^n$ be such that $t^n_{l} \to t_j$ and $t^n_{h} \to t$ (as usual we omit the dependence on $n$ in the indices $l$ and $h$). Then, taking the sum of \eqref{e.kanbis} for $k=h,\dots,l-1$, yields
\begin{equation} \label{diss}
   \F ( u^n_{l} , v^n_{l} ) \le \F ( u^n_{h} , v^n_{h} ) + \int_{t_{h}^n}^{t_{l}^n} \P^n (t, u^n(t) , v^n (t) ) \,dt + C \sum_{k=h}^{l-1} | t^n_{k+1} - t^n_k |^2 .  
\end{equation}
	
Let us pass the limit as $n \to +\infty$. By strong convergence, seen in the proof of Lemma \ref{ls.strong}, we infer that $\F ( u^n_{l} , v^n_{l} ) \to \F ( u(t_j) , v(t_j) )$ and $\F ( u^n_{h} , v^n_{h} ) \to \F ( u(t) , v(t) )$;  
indeed, we can write $v^n_{l}= ( v^n(t^n_{l}) - v^n(t_j) ) + v^n(t_j)$, where $v^n(t_j)\to v(t_j)$ in $H^1(\Omega)$ by Lemma \ref{l.strong}, and $\|v^n(t^n_{l})- v^n(t_j)\|_{H^1}\leq C|t_{l}^n-t_j|\to 0$ by Lipschitz continuity.  Similarly we get that   $u^n_l\to u(t_j)$ in $W^{1,p}(\Omega)$. We can apply the same reasoning to  $(u_h^n,v^n_h)$. 
If we take $j$ sufficiently large that $|t_j-t|<\eps$, then $ | t^n_{k+1} - t^n_k | < \eps$ and thus 
$$
	\sum_{k=h}^{l-1} | t^n_{k+1} - t^n_k |^2 \le \eps| t^n_{l} - t^n_{h}| .
$$
It now remains to prove that 
\begin{align*} 
\lim_{n\to +\infty}\int_{t^n_{h}}^{t^n_{l}} \P^{n}  ( t, u^n(t), v^n(t)) \, dt & = \int_{t}^{t_j }   \mathcal{P} (t, u(t), v(t)) \, dt .
\end{align*}
To this end, let us write
$$
	\P^{n} (t,u^n(t), v^n(t)) := \dot{a}^n (t) \int_\Omega \stress(u^n(t),v^n(t)) :  \strain( \hat g )  \, dx 
	:= \dot{a}^n(t)  P^n (t) . 
$$
We will now prove that $a^n  \weakstarto a:=\alpha\circ c$ in $W^{1,\infty }(0,T)$ and that $P^n \to P $ strongly in $L^1(0,T)$,  where $P(t):=\int_\Omega \stress(u(t),v(t)) :  \strain( \hat g )  \, dx$.
Let us start from $a^n$. It is enough to check the pointwise convergence of $a^n$ to $a$ in $[0,T]$. 
We start observing what happens on the points $t^n_k\in \Upsilon^n$:
\begin{equation}\label{tk}
  |a^n(t^n_k)-a(t^n_k)|=|\alpha(c^n(t^n_k))-\alpha(c(t^n_k))|\leq C |c^n(t^n_k)-c(t^n_k)| , 
\end{equation}
where we used the fact that $\alpha$ is Lipschitz continuous. 
Let us now see what happens on $I_u$. As we have proven in Lemma \ref{ccost}, $c$ is constant on the connected components of $I_u$, and so it is $a=\alpha\circ c$. For $t\in I_u$, there  exists $n$ sufficiently big such that $t\in (t^n_k, t^n_{k+1})\subset I_u $ (with $k$ depending on $n$). By definition of $a^n$, $a^n(t)=(1-\theta ^n)a^n(t^n_k)+\theta^na^n(t^n_{k+1})$ for some $\theta^n\in(0,1)$, while $a$, being constant on $[t^n_{k},t^n_{k+1}]$, can be written as any convex combination of $a(t^n_{k})$ and $a(t^n_{k+1})$. Hence:
$$|a^n(t)-a(t)|=(1-\theta^n)|a^n(t^n_k)-a(t^n_k)|+\theta^n|a^n(t^n_{k+1})-a(t^n_{k+1}))|\to 0,$$
following from \eqref{tk} and the fact that $c^n\to c$ uniformly in $[0,T]$. 
Lastly, for $t \in I_s$, we can take a sequence $\{t^n_{k}\in \Upsilon^n\}_{n\in \mathbb{N}}$ (with $k$ depending on $n$) converging to $t$ 
and therefore: 
\begin{align*}
    |a^n(t)-a(t)|&\leq |a^n(t)-a^n(t^n_k)|+  |a^n(t^n_k)-a(t^n_k)| +|a (t^n_k)-a(t)|
\\&\leq C |t-t^n_k| +  |a^n(t^n_k)-a(t^n_k)| \to 0, 
\end{align*}
where we used the fact that $a^n$ and $a$ are Lipschitz continuous together with \eqref{tk} and the fact that $c^n\to c$ uniformly in $[0,T]$.

Let us now check that $P^n \to P $ strongly in $L^1(0,T)$. This follows from the dominated convergence theorem: 
 $P^n  $ pointwise converges to $P $ 
 and it is uniformly bounded.
Indeed, we know that $v^n(t) \to v(t)$ strongly in $L^r(\Omega)$ and that $u^n (t) \to u (t)$ strongly in $W^{1,p}(\Omega)$ and therefore: 
	\begin{align*} 
\int_\Omega \stress(u^n(t),v^n(t)) :  \strain( \hat g )  \, dx  \to \int_\Omega \stress(u(t),v(t)) :\strain( \hat g)  \, dx. 
\end{align*} 
All in all, we get
$$
\F ( u(t_j) , v(t_j) ) \le \F ( u (t)  , v (t)  ) +  \int_{t}^{t_j} \int_\Omega  \stress(u(s),v(s)): \strain(\dot{a}(s)\hat g) \,dx\, ds   +\eps|t_j-t|  . 
$$ 
Dividing by $|t_j-t| $ and passing to the limit as $t_j \to t$ gives $\dot{\F} (t) \le \P (t, u(t) ,v(t)) + \eps $ and we conclude by arbitrariness of $\eps > 0$. \qed

\begin{corollary} \label{c.dernullz} $\partial_v \F ( u(t) ,  v (t) ) [ \dot{v}(t)] = 0$ a.e. in $I_s$.\end{corollary}
\proof This is a straightforward consequence of Theorem \ref{t.enbal}. Indeed:
$$  
        \dot{\F} (u(t),v(t)) = \partial_u \F ( u(t) ,  v (t) ) [ \dot{u} (t)]+\partial_v \F ( u(t) ,  v(t) ) [ \dot{v}(t)] = \partial_u \F ( u(t) ,  v (t) ) [ \dot{u} (t)] 
$$ 
a.e. in $I_s$ \bl and hence $\partial_v \F ( u(t) ,  v(t) ) [ \dot{v}(t)]=0.$ \qed

%% file: griffith.tex
\section{Griffith's criterion} \label{GR}

In this section we consider the elastic energy $W ( \strain , v ) = \psi (v) W_+ ( \strain ) + W_- ( \strain ) $ where 
$$ 
	W_+ (\strain) :=   \mu   | \strain_d |^2 + \kappa   | \strain_v^+ |^2 + \delta , 
	\quad
	W_- (\strain)  :=   \kappa   | \strain_v^- |^2  ,
$$
for some $\delta \ge 0$. Clearly, adding a constant to the energy does not change its mechanical properties. 
For $\delta > 0$ (arbitrarily small) a better description of the energy release holds (see Proposition \ref{p.non-deg}), while for $\delta =0$ we recover the classical phase-field elastic energy. 

\subsection{Energy release}
In this section we provide comments and properties on the notion of energy release employed in Theorem \ref{t.KKT}. 
In literature, the definition of {\it  phase-field 
energy release} $\mathcal{G}_{\text{\sl max}}( t ,v)$ is not univocal and in the following we will give it in the same vein as  \cite{KneesNegri_M3AS17}. Along the lines of the aforementioned article, let us define the reduced energy  
$$
	\tilde{\E} ( t ,v) := \E (  u_{t,v} , v)
	\quad \mbox{ for } \quad
	u_{t,v} \in \argmin \{ \E ( u , v) : u \in \U(t)  \} ,
$$
 that, as proven in the following Lemma, is differentiable as a function of $v$. \bl
\begin{lemma} \label{der.red}

The reduced energy functional  $\tilde{\E}$ is differentiable as a function of $v$ with respect to variations $\xi \in \Xi$.
Its derivative has the following integral representation:
\begin{equation} \label{gat}
    \partial_v\tilde{\E} (t, v)[\xi]=\partial_v \E (u_{t,v}, v)[\xi]=\int_{\Omega}v\xi \bigl( W_+(\strain(u_{t,v}))+ \delta \bigr)dx. 
\end{equation}
It follows that  $\partial_v \tilde{\E} (t,v) [\xi] \le 0$ for every $\xi \in 
\Xi$. 
\end{lemma} 
\proof
The directional derivative of $\tilde{\E} $  at $v$ is 
\begin{align*}
     \partial_v\tilde{\E} (t, v)[\xi]&=\lim_{h\to 0^+} \frac{\tilde{\E} (t, v+h\xi)-\tilde{\E} (t, v)}{h}=\lim_{h\to 0^+} \frac{\E (u_{t,v+h\xi}, v+h\xi)-\E (u_{t,v}, v)}{h}
\end{align*}
 Adding and subtracting the term $ \E(u_{t,v+h\xi},v)$, we obtain 
\begin{align*}
     &\lim_{h\to 0^+} \left( \frac{\E(u_{t,v+h\xi},v)-\E(u_{t,v},v)}{h}+\frac{\E(u_{t,v+h\xi},v+h\xi)-\E(u_{t,v+h\xi},v)}{h} \right). 
\end{align*} 
We examine the two addenda separately.
The first term is infinitesimal; indeed, using Corollary \ref{esp_quad} and Lemma \ref{le.KRZ2.5}, we obtain
\begin{align*} 
\frac{1}{h}\big|\E(u_{t,v+h\xi},v)-\E(u_{t,v},v)\big|\leq\frac{C}{h}\|u_{t,v+h\xi}-u_{t,v}\|^2_{H^1}\leq \frac{C}{h}\|h\xi\|^2_{L^r},
\end{align*}
that tends to zero as $h\to0$. As for the second term, we write it as follows:
\begin{align*}
     \frac{1}{2h}\int_{\Omega}  \bigl((v+h\xi)^2-v^2\bigr)  \big(W_+ (\strain(u_{t,v+h\xi}))+\delta\big)  \, dx &= \frac{1}{2h}\int_{\Omega}  \big(h^2\xi^2+2v\xi h \big) \big(W_+ (\strain(u_{t,v+h\xi}))+\delta\big)\, dx 
     \\&= \int_{\Omega}  \bigl( \tfrac12 h\xi^2+v\xi\bigr) \big(W_+ (\strain(u_{t,v+h\xi}))+\delta\big)\,  dx .
\end{align*}
From Lemma \ref{l.reg} we know that $u_{t,v+h\xi}$ is bounded in $W^{1,\tilde p}(\Omega,\mathbb{R}^2)$ uniformly with respect to $h$ and $\xi$. Choosing $1 \le q < +\infty$ such that $\tfrac{2}{q} + \tfrac{2}{\tilde p} =1$ we get 
$$\tfrac{1}{2}\int_{\Omega}  h\xi^2  \big(W_+ (\strain(u_{t,v+h\xi}))+\delta\big)  \, dx\leq
\tfrac12 h\big( ||\xi ||^2_{L^q}||u_{t,v+h\xi}||^2_{W^{1,\tilde p}}+\delta ||\xi ||^2_{L^2}\big)\leq  Ch\to 0. $$
Now consider $\int_{\Omega}  v\xi  \big(W_+ (\strain(u_{t,v+h\xi})) +\delta\big)\, dx $ and observe that  it converges to \eqref{gat}. 
By Lemma \ref{le.KRZ2.5}, indeed, 
\begin{align*}
    ||  v\xi W_+ (\strain(u_{t,v+h\xi})) & - v\xi W_+ (\strain(u_{t,v}))  ||_{L^1}=\int_{\Omega}  |v\xi   \big( W_+ (\strain(u_{t,v+h\xi})) - W_+ (\strain(u_{t,v})) \big)| \, \, dx \\&\leq C\|u_{t,v+h\xi} -u_{t,v}\|_{H^1}\to 0. 
\end{align*}
 The proof is concluded.  \qed

\noindent
We are now ready  to define the phase-field energy release rate for $v \in \V$ and $t \in [0,T]$ as:
\begin{align*}
	\mathcal{G}_{\text{\sl max}}(t,v) := \mathrm{sup} \big\{ \! - \partial_v \tilde{\E} (t, v) 
[\xi] : \xi \in \hat\Xi (v) 
\big\} , 
\end{align*} 
where  $\hat\Xi (v)$ is the set of {\it 
normalized admissible variations} with respect to the ''phase-field crack 
length'', i.e., 
$$
	\hat\Xi (v) = \{ \xi \in \Xi : d \mathcal{L} (v) [\xi]  = 1 \}  .
$$ 
First, note that $\mathcal{L} (v) = \tfrac12 \| v -1\|^2_{H^1}$ and thus $d \mathcal{L} (v) [\xi] = \langle v-1 , \xi \rangle_\V = 0$ for every $\xi$ when $v \equiv 1$. Therefore, we assume that $v \not \equiv 1$, so that 
$\hat\Xi(v)\neq\emptyset$. 
Clearly $\mathcal{G}_{\text{\sl max}}(t,v) \ge 0$ since $\partial_v \tilde{\E} (t,v) [\xi] \le 
0$ for every $\xi \in \hat\Xi (v)$.
\bl
To better understand this definition, it is important to observe that it is given in analogy with the sharp crack setting, where the maximal energy release is defined as the steepest descent of the (reduced) elastic energy with respect to crack elongation among any admissible direction. Indeed, if we consider the admissible variations $\xi\in \Xi$ such that $d \mathcal{L} ( v) [ \xi ] > 0$ and for each of them we define a sort of forward derivative of $\tilde{\E}$ with respect to the phase-field crack length, the following identities hold:
\begin{align}
   \sup \, \biggl\{ \! - \!\lim_{\,\,h\to 0^+} \frac{\tilde{\E}  (t,v+h\xi)-\tilde{\E} (t,v)  }{\mathcal{L} (v+h\xi)-\mathcal{L}(v)}  \,:\,\xi\in \Xi  \text{ with } d \mathcal{L} ( v) [ \xi] > 0  \biggr\} & \nonumber \\
   = \mathrm{sup} \,  \biggl\{ \! -  \frac{\partial_{v}\tilde{\E}(t,v)[\xi]}{d\mathcal{L} (v)[\xi]}  \,:\,\xi\in \Xi  \text{ with } d \mathcal{L} ( v) [ \xi] > 0   \biggr\} & \label{eo} \\
  =\mathrm{sup} \, \Bigl\{ \! -  \partial_{v}\tilde{\E} (t,v)[\xi] \,:\,\xi \in \Xi \text{ with } d \mathcal{L}  (v) [\xi]  = 1  \Bigr\} & = \mathcal{G}_{\text{\sl max}}( t, v) .  \nonumber
\end{align}
The second equality follows from linearity with respect to $\xi$, which allows to consider admissible variations normalized with respect to the phase-field crack elongation. 
In principle, given an arbitrary function $v \in \V$ it may happen that $d \mathcal{L} (v) [\xi] \le 0$ for some $\xi \in \Xi$; these variations would not be admissible in the representation of $\mathcal{G}$ given above. However, if $( u , v)$ is an equilibrium configuration of the energy (which is the case in quasi-static evolutions) we have $d \mathcal{L} (v) [\xi] \ge 0$ for every $\xi \in \Xi$ (independently of $\delta \ge 0$); moreover, if $\delta>0$ (but arbitrarily small) and if $( u , v)$ is an equilibrium configuration then $d \mathcal{L} ( v ) [ \xi ] > 0$ for every $\xi \in \Xi$ with $\xi \neq 0$, as stated in next Proposition.

\begin{proposition} \label{p.non-deg} Let $(u,v)$ be an equilibrium configuration for the energy. Then $d \mathcal{L} ( v ) [ \xi ] \ge 0$ for every $\xi \in \Xi$. Moreover, if $\delta >0$ then $d \mathcal{L} ( v ) [ \xi ] > 0$ for every $\xi \in \Xi$ with $\xi \neq 0$. 
\end{proposition}

\proof Since $\partial_v \F (u,v) [ \xi ] = \partial_v \E ( u , v ) [\xi] + G_c d \mathcal{L} ( v) [\xi] \ge 0$ in equilibrium configurations, it follows that $G_c  d \mathcal{L} ( v) [\xi] \ge - \partial_v \E ( u , v ) [\xi] \ge 0$ for every $\xi \in \Xi$.

Let $\delta >0$. Assume by contradiction that $d \mathcal{L} ( v ) [ \xi ] = 0$ for some $\xi \in \Xi$ with $\xi \neq 0$. It follows that $\partial_v \E ( u   , v   ) [ \xi ]  \ge 0$ and then $ \partial_v \E ( u   , v   ) [ \xi ]  = 0$. 
However, 
$$
	 \partial_v \E ( u   , v   ) [ \xi ]  = \int_\Omega v \xi ( W_+ ( \strain(u)) + \delta) \, dx = 0 
$$
implies that $v \xi =0$ a.e.~in $\Omega$ since $ W_+ ( \strain(u)) + \delta >0$. Since $v\in H^1(\Omega)$ we have $\nabla v = 0$ a.e.~in the set $\{ v = 0 \}$ (including the cases in which the measure of the set vanishes) and, similarly, $\nabla \xi = 0$ a.e.~in the set $\{ \xi = 0 \}$. As $v \xi =0$ a.e.~in $\Omega$ it follows that $\nabla v \cdot \nabla \xi = 0$ a.e.~in $\Omega$. 
As a consequence, 
\begin{align*}
	d \mathcal{L} ( v ) [ \xi ] & = \int_\Omega ( v-1) \xi + \nabla v \cdot \nabla \xi \, dx = \int_{ \Omega} - \xi \, dx = 0 ,
\end{align*}
which implies $\xi =0$ a.e.~in $\Omega$ since $\xi \le 0$.  \qed 

\noindent 
From the previous Lemma we get the following representation of the energy release. 

\begin{proposition} \label{Gequil} Let $\delta>0$ and let $(u,v)$ be an equilibrium configuration for the energy. Then 
\begin{align}
 \mathcal{G}_{\text{\sl max}}( t, v)  & = 
 	\mathrm{sup} \, \Bigl\{ \! -  \partial_{v}\tilde{\E} (t,v)[\xi] \,:\,\xi \in \Xi \text{ with } d \mathcal{L}  (v) [\xi]  = 1  \Bigr\} \nonumber  \\
 	&  = \mathrm{sup} \,  \biggl\{ \! -  \frac{\partial_{v}\tilde{\E}(t,v)[\xi]}{d\mathcal{L} (v)[\xi]}  \,:\, \xi\in \Xi,\, \xi \neq 0 \biggr\} & \label{equ-d} \\
 	& = \limsup_{z \to v^-}  \frac{ [ \tilde{\E}(t,z) - \tilde{\E}(t,v)]_- }{ \mathcal{L} (z) - \mathcal{L}(v) }, \label{equ-t}
 \end{align}
 where $z \to v^-$ means that $z \to v$ in $L^1(\Omega)$ with $0 \le z \le v$ and $z \neq v$. 
\end{proposition}

\proof The first identity follows immediately from the previous Proposition and \eqref{eo}. 
As $0 \le z \le v$ we have $\tilde{\E}(t,z) \le \tilde{\E}(t,v) $. 
Therefore
$$
\limsup_{z \to v^-}  \frac{ [ \tilde{\E}(t,z) - \tilde{\E}(t,v)]_- }{ \mathcal{L} (z) - \mathcal{L}(v) }
=
\limsup_{z \to v^-}  \frac{ \tilde{\E}(t,v) - \tilde{\E}(t,z)  }{ \mathcal{L} (z) - \mathcal{L}(v) } . 
$$
Note also that $\mathcal{L} (z) > \mathcal{L}(v)$ since $\mathcal{L}$ is quadratic and thus 
$$
	\mathcal{L} (z) = \mathcal{L} (v) + d \mathcal{L} (v) [ z - v] + \tfrac12 \| z - v \|^2_{H^1} ,  
$$
where $d \mathcal{L} ( v) [ z - v] >0$ by Proposition \ref{p.non-deg} (unless $z=v$).
\\
To check that 
$$
\mathrm{sup} \,  \biggl\{ \! -  \frac{\partial_{v}\tilde{\E}(t,v)[\xi]}{d\mathcal{L} (v)[\xi]}  \,:\, \xi\in \Xi,\, \xi \neq 0\bl\biggr\} \le \limsup_{z \to v^-}  \frac{ \tilde{\E}(t,v) - \tilde{\E}(t,z) }{ \mathcal{L} (z) - \mathcal{L}(v) }, 
$$it is not restrictive to assume that $\partial_{v} \tilde{\E}(t,v)[\xi] < 0$, otherwise $\partial_{v}\tilde{\E}(t,v)[\xi]=0$ and there is nothing to prove. Given $\xi \in \Xi$ let us choose $z = [ v + h \xi]_+$ and let $h \to 0^+$. First, note that $z \neq v$; indeed if $[ v + h \xi]_+ = v$ then $\xi$ would be supported in the set $\{ v = 0\}$, and thus $\partial_{v} \tilde{\E}(t,v)[\xi] = 0$. Then, $\mathcal{L} (v ) < \mathcal{L} (z) \le \mathcal{L} ( v + h \xi)$ and thus
$$
	\frac{1}{\mathcal{L} ( v + h \xi) - \mathcal{L}(v)} \le \frac{1}{\mathcal{L} (z) - \mathcal{L} ( v )} .
$$
Moreover, $(v+h\xi)^2 \ge z^2$ and thus $\tilde{\E} ( t , v + h \xi) \ge \tilde{\E}(t,z)$. Hence,
$$
	\tilde{\E} ( t , v ) - \tilde{\E}(t, v + h \xi) \le \tilde{\E} ( t , v ) - \tilde{\E}(t, z)  = [ \tilde{\E} ( t , z ) - \tilde{\E}(t, v) ]_- . 
$$
In conclusion,
$$
	\frac{\tilde{\E} ( t , v ) - \tilde{\E}(t, v + h \xi) }{\mathcal{L} ( v + h \xi) - \mathcal{L}(v)} \le \frac{[ \tilde{\E} ( t , z ) - \tilde{\E}(t, v) ]_- }{\mathcal{L} (z) - \mathcal{L} ( v )} .
$$
Passing to the limsup as $h \to 0^+$ gives the required inequality, by arbitrariness of $\xi$. 

It remains to prove the opposite inequality. To this end, it is not restrictive to assume that 
$$
	L = \limsup_{z \to v^-}  \frac{ \tilde{\E}(t,v) - \tilde{\E}(t,z) }{ \mathcal{L} (z) - \mathcal{L}(v) } > 0 . 
$$
First, consider the case $L < +\infty$. For $\eps>0$ (sufficiently small) there exists $z_n \to v^-$ such that 
$$
	\frac{ \tilde{\E}(t,v) - \tilde{\E}(t,z_n)  }{ \mathcal{L} (z_n) - \mathcal{L}(v) } \geq L - \eps > 0.
$$
Note that $ \tilde{\E} (t, z_n) \to \tilde{\E}(t, v)$, and thus $ \mathcal{L} (z_n) \to \mathcal{L}(v)$ (otherwise the limit of the left-hand side would vanish). As a consequence $z_n \to v$ strongly in $H^1(\Omega)$ and thus in $L^q(\Omega)$ for any $1 \le q < +\infty$. 
Setting $u_n \in \argmin \{  \E  (  u , z_n ) : u \in \U (t)\}$ and $u \in \argmin \{  \E (  u , v ) : u \in \U (t)\}$, let us write
\begin{align*}
	\tilde{\E} (t, z_n) 
		& = \E ( u_n, z_n ) - \E ( u, z_n ) + \E (u, z_n) \geq \E ( u, z_n ) -C \| u_n - u \|_{H^1}^2 \\&\ge \E ( u, z_n ) - C\| z_n - v \|^2_{L^r},
\end{align*}
where once again we used Lemma \ref{esp_quad} and Lemma \ref{le.KRZ2.5}. By convexity
$$ \E ( u, z_n ) 
		\ge \E ( u, v) + \partial_v \E ( u, v ) [z_n - v] 
		= \tilde{\E} ( t, v) +  \partial_v \tilde{\E} ( t, v )  [z_n - v] . 
$$
Hence, 
\begin{align*}	
	\tilde{\E} (t, v) - \tilde{\E} (t, z_n) \le - \partial_v \tilde{\E} ( t , v ) [z_n - v]  +  C \| z_n - v \|^2_{ L^r} . 
\end{align*}
Since $\mathcal{L}$ is quadratic we can write
$$
	\mathcal{L} (z_n) - \mathcal{L} (v) = d \mathcal{L} (v) [ z_n - v ] + \tfrac12 \| z_n - v \|_{H^1}^2 
$$
and then
\begin{equation} \label{e.eleps}
	L - \eps \le \frac{- \partial_v \tilde{\E} ( t , v ) [z_n - v]  +  C \| z_n - v \|^2_{L^r} }{d \mathcal{L} (v) [ z_n - v ] + \tfrac12 \| z_n - v \|_{H^1}^2}  . 
\end{equation} 
Note that $- \partial_v \tilde{\E} ( t , v ) [z_n - v] \ge 0$ and let $\xi_n = ( z_n -v) / \| z_n - v \|_{L^r}$, so that $\| \xi_n \|_{L^r} = 1$. 
\begin{itemize}
    \item If $\| \xi_n \|_{H^1} \le C$ then (up to non-relabelled subsequences) $\xi_n \weakto \xi$ in $H^1(\Omega)$ and thus $\xi_n \to \xi$ in $L^r(\Omega)$ by Sobolev embedding. In particular $\xi \neq 0$ since $\| \xi \|_{L^r} =1$. 
From \eqref{e.eleps} we get 
\begin{align*}
	L - \eps 
		& \le  \frac{- \partial_v \tilde{\E} ( t , v ) [ z_n - v ]}{d \mathcal{L} (v) [ z_n - v ] }  + C \frac{ \| z_n - v \|^2_{L^r} }{d \mathcal{L} (v) [ z_n - v  ]} \\
		& \le \mathrm{sup} \,  \biggl\{ \! -  \frac{\partial_{v}\tilde{\E}(t,v)[\xi]}{d\mathcal{L} (v)[\xi]}  \,:\, \xi\in \Xi,\, \xi \neq 0 \biggr\} + C \frac{ \| z_n - v \|_{L^r}}{d \mathcal{L} (v) [ \xi_n ]}.
\end{align*}
In the limit as $n \to +\infty$ the last term vanishes since $d \mathcal{L} (v) [ \xi_n ] \to d \mathcal{L} (v) [ \xi ] > 0 $, by Proposition \ref{p.non-deg}, while $z_n \to v$ in $L^r (\Omega)$.
\item If $\| \xi_n \|_{H^1} \to +\infty$ for a (non-relabelled) subsequence then \eqref{e.eleps} yields
\begin{align*}
	L - \eps & \le  \frac{- \partial_v \tilde{\E} ( t , v ) [ z_n - v ]}{d \mathcal{L} (v) [ z_n - v ] }  + C \frac{ \| z_n - v \|^2_{L^r} }{\| z_n - v \|_{H^1}^2} \\
		& \le \mathrm{sup} \,  \biggl\{ \! -  \frac{\partial_{v}\tilde{\E}(t,v)[\xi]}{d\mathcal{L} (v)[\xi]}  \,:\, \xi\in \Xi,\, \xi \neq 0 \biggr\} + \frac{  C  }{\| \xi_n \|_{H^1}^2} .
\end{align*} 
Clearly last term vanishes in  the limit as $n \to +\infty$. In both cases we conclude by the arbitrariness on $\eps$.
\end{itemize}
In the case $L = +\infty$ it is enough to replace $L-\eps$ with an arbitrary large value.  \qed

\begin{remark} \normalfont Note that writing the energy release as the ``slope''
$$
	\mathcal{G}_{\text{\sl max}}(t, v) =  \limsup_{z \to v^-}  \frac{ [ \tilde{\E}(t,z) - \tilde{\E}(t,v)]_- }{ \mathcal{L} (z) - \mathcal{L}(v) },
$$
gives the most general way of looking at variations of energy with respect to variations of phase-field crack length. In perspective, this would be the natural definition of energy release for the $\Gamma$-limit functional in the space $SBD$, where directional derivatives seems not general enough to describe all the possible unilateral variations of the crack.
\end{remark}

\subsection{Griffith's criterion} 

Before presenting the phase-field version of Griffith's criterion satisfied by the limit evolution $v$, given by Lemma \ref{l.strong}, we shall study the energy functional on the set $I_u$.

\begin{lemma} \label{f10} Let $(a,b)$ be a connected component of $I_u$, then $\F ( u (b) , v(b) ) \le \F ( u(a) , v(a))$. Moreover 
$$	\dot{\F} ( u(a^+) , v(a^+) ) \ge 0  \,, \quad \dot{\F} ( u(b^-) , v(b^-) ) = 0  .    $$
\end{lemma}

\proof 
We recall that $v$ is affine in $(a,b)$ and denote for convenience
\begin{equation} \label{f}
f(t):=\F(u (t) , v (t) ),  \quad 
e (t) := \tilde\E ( t , v(t) ) = \E ( u(t) , v(t)) , \quad \ell(t) := \mathcal{L} ( v(t) ) . 
\end{equation}
In this notation, we will show that the function $f$ satisfies $f(b) \leq f(a)$,  $\dot{f} (a^+) \ge 0$ and $\dot{f}(b^-)=0$.

As we have already mentioned, since $a,\,b \in I_s$ by strong convergence (seen in the proof of Lemma \ref{ls.strong}) combined with Lemma \ref{conv_ajbj} 

there exists $u^n(t^n_{k+1})\to u(b)$ (as usual $k = k(n)$) and $u^n(t^n_k)\to u(a)$ strongly in $W^{1,p} (\Omega ; \mathbb{R}^2)$, $v^n(t^n_k)\to v(a)$ and $v^n(t^n_{k+1})\to v(b)$ strongly in $H^1 (\Omega)$. 
It follows that $\F(u(t^n_{k}),v(t^n_{k}))\to \F(u(a),v(a))=f(a)$ and  $\F(u(t^n_{k+1}),v(t^n_{k+1}))\to \F(u(b),v(b))=f(b)$.
By minimality together with \eqref{Pnk}
\begin{align}  \label{e.64}
     \F (u^n_{k+1} , v^n_{k+1}) & \le   \F(u^n_{k}, v^n_{k})+\int_{t_{k}^n}^{t_{k+1}^n}\P^{n} (t,u^n_k , v^n_k ) \,dt + C | t^n_{k+1} - t^n_k |^2 , 
\end{align}
where  
$$	\P^{n} (t,u, v) := \int_\Omega \stress(u,v) :  \strain( \dot a^n(t) \hat g )  \, dx . $$ 
Now, by definition of $a^n$, we have that, for $t\in (t_{  k},t_{  k+1})$, $\dot{a}^n (t)= \frac{\alpha(s_{  k+1}^n)-\alpha( s_{  k}^n)} {\tau_{k+1}}$ and thus the integral of the power  in \eqref{e.64} reads 
$$ \int_\Omega \big(   \alpha(s_{  k+1}^n)-\alpha( s_{  k}^n) \big) \stress(u^n_k,v^n_k):\strain(\hat{g}) \, dx \, \to 0$$
 since $\alpha$ is continuous and $s_{  k+1}^n=s_{  k}^n+\tfrac{1}{n}\to s_{  k}^n$. Therefore  
$$f(b)  = \lim_{n \to +\infty}   \F (u^n_{k+1} , v^n_{k+1}) \leq \lim_{n\to +\infty}\F(u^n_{  k}, v^n_{  k}) =f(a).$$

Let's now prove that $\dot{f} (a^+)=0$. We know  that $v^n_{  k}\to v(a )$ and $v^n_{  k+1} \to v(b )$  strongly in $H^1 (\Omega)$. 
Hence
$$\dot{v}^n_{  k+1}=\frac{v^n_{  k+1}-v^n_{  k}}{t^n_{  k+1}-t^n_{  k}}\longrightarrow\frac{v(b )-v(a )}{b -a }=
\dot{v} (t) \quad  \text{strongly in $H^1(\Omega)$ for all $t \in(a,b)$.}$$ 
Being $( u(a), b(a) )$ an equilibrium configuration (see Theorem \ref{EQ}) and being $\dot{v} (a^+) \le 0$ an admissible variation, \eqref{equ} yields $\dot{f} (a) \ge 0$. 

It remains to prove that $\dot{f} (b^-)=0$. 
Recall that, by \eqref{eq0DISC}, $\partial_v\F(u^n_{k+1} ,v^n_{k+1} )[\dot{v}^n_{k+1} ]=0$.
Its explicit form is:
$$\partial_v\F(u^n_{  k+1} ,v^n_{  k+1} )[\dot{v}^n_{  k+1} ]=\int_{\Omega}v^n_{  k+1} \dot{v}^n_{  k+1}  \bigl( W_+(\strain(u^n_{  k+1})+\delta \bigr) \, dx+G_c\int_{\Omega}(v^n_{  k+1} -1)\dot{v}^n_{  k+1} +\nabla v^n_{  k+1} \nabla \dot{v}^n_{  k+1} = 0 . $$
The strong convergence of $v^n_{k+1}$ and $\dot{v}^n_{k+1}$ combined to the fact that $u^n_{  k+1}\to u(b )$ strongly in $W^{1, p}(\Omega)$, gives 
$$0=\partial_v\F(u^n_{  k+1} ,v^n_{  k+1} )[\dot{v}^n_{  k+1} ] \to \partial_v\F(u (b) ,v (b) )[\dot{v} (b^-)] = \dot{f}(b^-) ,$$ 
which concludes the proof. \qed

As we anticipated, neither the energy release nor the thermodynamic consistency enter into the definition of the time discrete evolution \eqref{e.newsequence}, however, the limit (time continuous) evolution satisfies both Griffith's criterion, in terms of the phase-field energy release, and thermodynamic consistency, as stated in Theorem \ref{t.KKT} (recalled hereafter for reader's convenience). 

\medskip

{\it Assume that $v_0 \not \equiv 1$. 
The limit evolution $v$, obtained in the limit of the time discrete evolution \eqref{e.newsequence}, satisfies $0 \le v (t)  \le 1$ for every $t \in [0,T]$, the irreversibility constraint $\dot{v} (t) \le 0$ and the thermodynamic consistency condition $\dot{\ell} ( t ):=d\mathcal{L}(v(t))[\dot{v}(t)] \ge 0$ a.e.~in $[0,T]$. 

Moreover, an extended version of Griffith's criterion holds: 
\begin{itemize}
\item $\mathcal{G}_{\text{\sl max}}( t  , v(t) ) \le G_c$ everywhere in $I_s$, 
\item $(\mathcal{G}_{\text{\sl max}}( t , v(t) ) - G_c ) \, \dot{\ell} (t) = 0$ a.e.~in $I_s$,
\item  $\partial_v \tilde\E ( t , v )  [ \dot{v}(t)]  = -  \dot{\ell}(t) \, \mathcal{G}_{\text{\sl max}} ( t ,v(t) )$ a.e.~in $I_s$,
\item $\mathcal{G}_{\text{\sl max}}(t, v(t)) \ge G_c$ on a subset of positive measure of $I_u$.

\end{itemize}
}
\bl 
\noindent
{\bf Proof of Theorem \ref{t.KKT}.} {\bf I.} The validity of the  thermodynamic consistency condition on $I_s$ is a straightforward consequence of  Corollary \ref{c.dernullz}. Indeed, by definition $\partial_v \F ( u(t) ,  v (t) ) [ \dot{v}(t)] =\partial_v \E ( u(t) ,  v (t) ) [ \dot{v}(t)]+G_c \, d \mathcal{L} ( 
v(t) ) [\dot{v}(t)]$. Now, from the aforementioned Corollary, we know that $\partial_v \F ( u(t) ,  v (t) ) [ \dot{v}(t)] = 0$ a.e. in $I_s$
and therefore:
$$
	 G_c \, d \mathcal{L} ( 
v(t) ) [\dot{v}(t)]= -\partial_v \E ( u(t) , v(t)) [\dot{v}(t)] =-\int_\Omega v(t) \,  \dot{v}(t)\, W_+(\strain(u(t) )  \, dx \ge 0 .  
$$ 
\\To prove the same property on $I_u$, we proceed as follows.  Let $(a,b)$ be a connected component of $I_u$. 
The function $t \mapsto e (t) = \tilde{\E} (t,v(t)) = \E ( u(t) , v(t)) $ is non increasing in $(a,b)$. Indeed, if $t_1< t_2$ then $v ( t_1 ) \geq v ( t_2)$ and hence:
\begin{align*}  \E ( u (t_1) , v ( t_1) ) & = \int_\Omega \big( v^2  (t_1) +\eta \big) W_+(\strain(u (t_1) )  \, dx 
\geq \int_\Omega \big( v^2  (t_2) +\eta \big) W_+(\strain(u (t_1) ) \, dx  \\
& \geq \int_\Omega \big( v^2 (t_2) +\eta \big) W_+(\strain(u (t_2) ) \, dx  = \E ( u (t_2) , v  (t_2) ),
\end{align*} 
where the last inequality follows from the minimality of $u (t_2)$ (we recall that in $(a,b)$ the boundary condition is fixed since the control $c$ is constant, see Lemma \ref{ccost}). 
Now, writing $f(t) = e(t) + G_c \ell(t)$ as in \eqref{f},  from Lemma \ref{f10} we get 
$$ 0= \dot{f} (b^-)= \dot{e} (b^-) + G_c \, \dot{\ell} ( b^- )  .$$
Therefore, since the energy $t \mapsto e(t)$ is non increasing  $ d \mathcal{L}(v(b))[\dot{v}(b^-)] = \dot{\ell} (b^-) = - \dot{e} (b^-) / G_c  \geq 0.$
Similarly,  from the equilibrium condition in $a$ 
$$	
	0\leq \dot{f} (a ^+)=  \dot{e} (a^+) + G_c \, \dot{\ell} ( a^+ ) 
$$
and from the fact that $\dot{e} (a^+)\leq 0$, it follows that $d \mathcal{L}(v(a))[\dot{v}(a^+)]= \dot{\ell} ( a^+ ) \ge - \dot{e} (a^+) / G_c \geq0.$
Now, since $\ell$ is quadratic, its derivative is linear and since $\dot{\ell} (t) \ge 0 $ for both $t=a$ and $t=b$, it must be so also for every $t \in (a,b)$. We thus have that $\dot{\ell} (t) = \partial_v\mathcal{L}(v(t))[\dot v(t)]\geq0$ also for all $t \in I_u$. 

\medskip
\noindent {\bf II.} We will now prove that a classical version of Griffith's criterion holds on $I_s$, following \cite{KneesNegri_M3AS17}. 

\begin{itemize}
  \item  To prove that $\mathcal{G}_{\text{\sl max}}( t , v(t) ) \le G_c $ in $I_s$, 
we use the fact that, by \eqref{equ},
$$
	\partial_v \F ( u(t) , v(t)) [\xi] =\partial_v \E ( u(t) , v(t)) [\xi]  + G_c \, d \mathcal{L} ( v(t)) [\xi] \ge 0 \quad \mbox{ for every $\xi \in \Xi$}.  $$
In particular, if we take $\hat\xi \in \hat \Xi (v(t))$ we get $G_c  \ge  - \partial_v \tilde\E ( t,  v(t)) [\hat\xi]  
$
and thus 
\begin{equation*}
	G_c \ge \sup \big\{\!-\partial_v \tilde\E ( t, v(t)) [\hat\xi] : 
\hat\xi \in \hat\Xi (v(t))\big\} = \mathcal{G}_{\text{\sl max}}( t , v(t)) . 
\end{equation*}	
\item 
Finally, let us prove that $(\mathcal{G}_{\text{\sl max}}( t , v(t) ) - G_c ) \, \dot{\ell}(t) = 0$ a.e.~in $I_s$. Clearly, if $\dot{\ell}(t) =0$ there is nothing to prove. Consider $ \dot{\ell} (t) = d \mathcal{L} ( v(t) ) [ \dot{v}(t) ] >0$. Again by Corollary \ref{c.dernullz} we have that  
$\partial_v \F ( u(t) ,  v (t) ) [\lambda \dot{v}(t)] = 0$ for every  $\lambda > 0$ a.e.~in $I_s$. 
On the other hand, 
$
	\partial_v \F ( u(t) , v(t) )  [ \xi ]  \ge 0 \, \mbox{ for 
every $\xi \in \Xi$} 
$
and thus for every $\lambda > 0$
\begin{align*}
      \lambda \dot{v}(t) & \in  \argmin \{ 
      \partial_v \F (  u(t) , v(t) )  [\xi ]  :  \xi \in \Xi \}  \\  
	 & \in \argmin \{    \partial_v \tilde\E ( t , v(t))  [ \xi]  
	 + G_c d \mathcal{L} ( v(t) ) [\xi]  :  \xi \in \Xi \} . 
\end{align*}
Choosing $\lambda =  1 / d \mathcal{L} ( v) [ \dot{v}(t) ] = 1 / \dot\ell(t)$, in such a way that $ d \mathcal{L} ( v) [\lambda \dot{v}(t) ] =1$ (so that $\lambda \dot{v}(t)\in \hat\Xi(v(t))\subset\Xi$), we get that
\begin{align*}
      \lambda \dot{v}(t) &\in \argmin \{    \partial_v \tilde\E ( t , v(t) )  [ 
\hat\xi]+ G_c\, d \mathcal{L} ( v(t) ) [\hat\xi]   :  \hat\xi \in \hat\Xi(v(t)) \}  \\
 &\in \argmin \{    \partial_v \tilde\E ( t , v(t) )  [ 
\hat\xi]   :  \hat\xi \in \hat\Xi(v(t)) \} ,
\end{align*}since $ G_c \,d \mathcal{L} ( v(t) ) [\hat\xi]= G_c$ is constant in $\hat\Xi(v(t)) $. Therefore, by definition $$ \partial_v \tilde\E ( t , v )  [ \lambda \dot{v}(t)] 
= - \mathcal{G}_{\text{\sl max}}( t ,v(t) )$$ and
$$
	 \partial_v \F ( u(t) , v(t) )  [  \lambda \dot{v}(t) ] 
	 =  - \mathcal{G}_{\text{\sl max}}( t  , v(t)) + G_c = 0 ,
$$
which concludes the proof.
\end{itemize}

\noindent {\bf III.}  Let's now focus on the set of instability points $I_u$. Let $(a,b)$ be a connected component of $I_u$. By Lemma \ref{f10} we have $f(a) \ge f(b)$, therefore the subset of $(a,b)$ where $\dot{f} \le0$ is of positive measure. 

First, note the following: if the evolution is not identically constant (in which case there is nothing to prove) then $\dot{\ell} (t) > 0$ in $I_u$. Indeed, being $\dot{\ell} (t)$ linear and non-negative, if $\dot{\ell} (t)=0$ for some $t \in (a,b)$, then $\dot{\ell}$ should be zero everywhere and $\ell$ should be constant on $(a ,b )$. This means that $v$ would belong to a contour line of $\mathcal{L}(v)=
\tfrac12 \| v -1\|^2_{H^1(\Omega)}$. 
Now, by uniform convexity of $H^1 (\Omega)$, if $v (a) \neq v (b)$, then for any $t \in (a,b)$ we would have $\| v (t)  - 1  \|^2_{H^1}<\|v (a) -1 \|^2_{H^1}=\| v (b) - 1 \|^2_{H^1}$, i.e., $v (t)$ could not belong to the contour line of $\mathcal{L}$, that is absurd.

We will show that $\mathcal{G}(t,v(t))\geq G_c$ where $\dot{f} (t) \leq 0$. Write 
$$
 \dot{f} (t)=\dot{e} (t) +G_c \, \dot{\ell}(t) \leq 0.
$$
As we have seen, $ \lambda = \dot{\ell} (t) >0$, so we can divide, obtaining
\begin{equation}
G_c=G_c \frac{\dot{\ell} (t) }{\lambda}\leq - \frac{\dot{e} (t)}{\lambda}  = 
- \partial_v \tilde \E(t,v(t)) \Big[ \frac{\dot{v}(t)}{\lambda} \Big]  \leq \mathcal{G}_{\text{\sl max}}(t,v(t)) ,
\end{equation}
where in the last inequality we used the defintion of $\mathcal{G}_{\text{\sl max}}$ and the fact that $\dot{v} (t) / \lambda \in \hat\Xi (v (t))$, since by construction $\dot{\ell} (t) = d \mathcal {L} ( v(t)  ) [ \dot{v} (t) ] = \lambda$.
\qed

%% file: pde.tex

\bl 
\section{System of PDEs}

\begin{lemma} \label{l.measure} Let $\Omega \subset \mathbb{R}^2$ be a bounded, Lipschitz domain. Let $\zeta \in (H^1(\Omega))^*$ such that $\langle \zeta , \xi \rangle  \ge 0$ for every $\xi \in H^1( \Omega)$ with $\xi \le 0$. Then $\zeta$ is ``represented'' by $\mu \in \mathcal{M}^- (\bar\Omega)$ (the set of finite, negative, Radon measures supported in $\bar\Omega$). More precisely, 
$$
	\langle \zeta , \xi \rangle = \int_{\Omega} \xi d \mu_{| \Omega} + \int_{\partial \Omega} \xi d \mu_{| \partial \Omega} 
	\quad
	\text{for every $\xi \in H^1 ( \Omega) \cap C (\bar\Omega)$.}
$$
\end{lemma}

\proof Let $\tilde\Omega$ be an open set with $\bar{\Omega} \subset \tilde{\Omega}$ and consider the extension $\tilde\zeta\in  H^{-1} (\tilde{\Omega})$ given by $\langle \tilde{\zeta} , \phi \rangle = \langle \zeta , \phi_{| \Omega} \rangle$. Then $\langle \tilde\zeta , \phi \rangle \ge 0$ for every $\phi \in C^\infty_c (\tilde\Omega)$ with $\phi \le0$, since $\xi = \phi_{|\Omega} \in H^1(\Omega)$ and $\xi \le 0$.  As a consequence $\tilde\zeta$ (as a distribution) is represented by $\tilde\mu \in \mathcal{M}_{\mathrm{loc}}^- ( \tilde\Omega)$. Clearly, the support of $\tilde\mu$ is contained in $\bar\Omega$ and thus, being $\tilde{\mu}$ locally finite, it turns out that $\tilde{\mu} = \tilde{\mu}_{| \bar\Omega}$ is finite. 

Let $\xi \mapsto \tilde\xi$ be a (linear and continuous) extension from $H^1(\Omega)$ to $H^1_0 (\tilde\Omega)$. Then $\langle \zeta , \xi \rangle =  \langle \zeta , \tilde\xi_{|\Omega}  \rangle = \langle \tilde\zeta , \tilde\xi \rangle$. In particular, the integral representation holds for $\xi \in H^1 ( \Omega) \cap C (\bar\Omega)$.   \qed 

If $z \in H^1 (\Omega)$ with $\Delta z(t) \in L^2(\Omega)$ then $\partial z / \partial n$ (in weak form) denotes the operator  given by Green's formula 
$$
	\langle \partial z / \partial n , w \rangle  = \int_\Omega w \,  \Delta z  + \nabla w \cdot \nabla z \, dx 
$$
for $w \in H^1(\Omega)$. In a similar way,  if $z \in H^1(\Omega)$ such that $\Delta z$ (as a distribution) is represented by a finite measure, then we define the linear operator 
\begin{equation} \label{e.dudn}
	\langle \partial z / \partial n , w \rangle = \int_\Omega w \,  d (\Delta z)  + \int_\Omega \nabla w \cdot \nabla z \, dx ,
\end{equation} 
for $w \in H^1 ( \Omega) \cap C (\bar\Omega)$. 
Note that the first integral is well defined since $w \in C (\bar\Omega)$. 

\medskip
Let us denote by $\zeta (t) \in (H^1(\Omega))^*$ the linear operator $\langle \zeta (t) , \xi \rangle = \partial_v F (t, u(t) , v(t) ) [ \xi ]$. For $t\in I_s$ we have that $\langle \zeta (t) , \xi \rangle \ge 0$ for every $\xi \in H^1(\Omega)$ with $\xi \le 0$ and thus, by Lemma \ref{l.measure}, $\zeta (t) $ is represented by a measure $\mu (t) \in \mathcal{M}^- ( \bar\Omega)$. Let us write $\mu (t)  = \mu_{|\Omega} (t) +\mu_{|\partial\Omega} (t)$.

\begin{lemma} \label{l.measures}
For $t\in I_s$ we have $\Delta v(t) \in \mathcal{M} (\Omega)$ (the space of finite Radon measures in $\Omega$). Moreover, 
\begin{gather*}
	\mu_{|\Omega} (t)  =  - G_c \, \Delta v(t) + 2 v(t) W_+ \big( \strain ( u (t) ) \big) \,  + G_c (v(t)-1) \le 0 , \\
	\mu_{| \partial \Omega} (t) = G_c \, \partial v (t) / \partial n \le 0 . 
\end{gather*}
\end{lemma}

\proof
For $\phi \in C^\infty_c (\Omega)$ we have 
\begin{align*}
	\int_\Omega \phi \, d \mu_{|\Omega} (t) & = \langle \zeta (t)  , \phi \rangle = \partial_v F (t, u(t) , v(t) ) [ \phi ] \nonumber \\ 
	& = \int_\Omega 2 v(t) \phi \, W_+ \big( \strain ( u (t) ) \big)  \, dx + G_c \int_\Omega \nabla v(t)  \cdot \nabla \phi + (v(t)-1) \phi  \, dx \nonumber \\ 
	& = G_c \phantom{i}_{\mathcal{D}'} \langle - \Delta v(t) , \phi \rangle_{\mathcal{D}} + \int_\Omega [ 2 v(t) W_+ \big( \strain ( u (t) ) \big) \,  + G_c (v(t)-1) ] \phi  \, dx . 
\end{align*}
Hence
\begin{equation} \label{e.measures} 
	 G_c \phantom{i}_{\mathcal{D}'} \langle \Delta v(t) , \phi \rangle_{\mathcal{D}} = 
	 \int_\Omega [ 2 v(t) W_+ \big( \strain ( u (t) ) \big)  \,  + G_c  (v(t)-1) ] \phi  \, dx  - \int_\Omega \phi \, d \mu_{|\Omega} (t) .
\end{equation}
It follows that  $\Delta v(t) \in \mathcal{M} (\Omega)$ and that 
\begin{equation} \label{e.muomega}
	\mu_{|\Omega} (t) = - G_c \Delta v(t) + 2 v(t) W_+ \big( \strain ( u (t) ) \big)  \,  + G_c (v(t)-1) . 
\end{equation} 
Hence, for every $\xi \in H^1 (\Omega) \cap C (\bar\Omega)$ we can write 
\begin{align*}
	\langle \zeta (t) , \xi \rangle & =  \partial_v F (t, u(t) , v(t) ) [ \xi ] 
		= \int_\Omega G_c \nabla v(t) \cdot \nabla \xi + 2 v (t) W_+ \big( \strain( u(t)) \big) \xi + G_c ( v(t) - 1) \xi \, dx  \\
		& = G_c  \langle \partial v (t)  / \partial n , \xi \rangle - G_c \int_\Omega \xi \, d ( \Delta v(t) ) + \int_\Omega 2 v (t) W_+ \big( \strain( u(t)) \big) \xi + G_c ( v(t) - 1) \xi \, dx \\
		& = G_c \langle \partial v (t)  / \partial n , \xi \rangle + \int_\Omega \xi \, d \mu_{|\Omega} (t) . 
\end{align*}
Since
$$
	\langle \zeta (t) , \xi \rangle = \int_{\partial \Omega} \xi \, d \mu_{|\partial \Omega} (t)  + \int_\Omega \xi \, d \mu_{|\Omega} (t) , 
$$
it follows that $G_c\, \partial v(t) / \partial n$ is represented by the negative measure $\mu_{|\partial \Omega} (t)$. 
\qed

\begin{proposition} \label{p.KT2} Let $t \in I_s$ such that $\partial_v F ( t, v(t) , u(t) ) [ \dot{v}(t) ] = 0$.  If $\dot{v} (t) \in H^1 (\Omega) \cap C ( \bar \Omega)$ then 
\begin{equation} \label{e.KT2}
	 \big( - G_c \Delta v(t) +  2 v(t) W ( D u (t)) \,  + G_c (v(t)-1) \big)  \dot{v} (t) = 0 , 
	 \qquad
	 ( \partial v(t) / \partial n ) \dot{v}(t) = 0  , 
\end{equation}
in the sense of measures. 
\end{proposition}

\proof By Lemma \ref{l.measure} 
$$
	\partial_v F ( t, v(t) , u(t) ) [ \dot{v}(t) ] = \int_{\Omega} \dot{v}(t) d \mu_{| \Omega} (t)  + \int_{\partial \Omega} \dot{v}(t) d \mu_{| \partial \Omega} (t)  = 0 . 
$$
Since $\mu_{| \Omega} (t) \le 0$, $\mu_{| \partial \Omega} (t) \le 0$ and $\dot{v}(t)\le 0$ it follows that both the above integrals are non-negative, and thus they both vanish. By Proposition \ref{l.measures} we finally get \eqref{e.KT2}. \qed 

 Proposition \ref{l.measures} and \ref{p.KT2} give Theorem \ref{t.PDE}.

%% file: alternate.tex

\section{Alternate minimization} 
\label{M-ele}

Given $v_0 \in \V$,  with $0  \le v_0 \le 1$ let $u_0 \in \argmin \{  \F  (  u , v_0) : u \in \ \U (0) \} $. For $n\in \N$ set $s^n_k := k / n $ for $k=0,...,n$ and $v^n_0 : = v_0$. In this section we consider 
the alternate minimization scheme \cite{BourdFrancMar00} at $s^n_{k}$ with initial condition $v^n_{k-1}$  (defined in $s^n_{k-1}$ and obtained at the previous iteration), i.e., set $v^n_{k,0} := v^n_{k-1}$ and $u^n_{k,0} := u^n_{k-1}$, then we define by induction for $i\geq 0$: 
\begin{equation}  \label{e.altminM}
\begin{cases}
	u^n_{k,i+1} \in \argmin \{  \F  (  u , v^n_{k,i} ) : u \in \U (s_k^n)
\} 
\\
         v^n_{k,i+1} \in \argmin \{  \F ( u^n_{k,i+1} , v ) : v \in \V 
\mbox{ with } v \le v^n_{k,0} = v^n_{k-1} \} . 
\end{cases} 
\end{equation} 
The updates $u^n_{k}$ and  $v^n_{k}$ are respectively defined as $u^n_k:= \lim_{i \to \infty} u^n_{k,i}$ and $v^n_k:= \lim_{i \to \infty} v^n_{k,i}$; existence of the limits (up to non-relabelled subsequences) is proved in the following Proposition. 

\begin{remark} \normalfont As in the original and mostly used version of the alternate minimization algorithm, we require that $v^n_{k+1} \leq v_k^n$ for every $k=0,\dots n-1$,\, and not $v^n_{k,i+1} \leq v_{k,i}^n$ for every $i \in \mathbb{N}$ (i.e., monotonicity over each staggered iteration)  as in \cite{KneesNegri_M3AS17}.
\end{remark}

\begin{proposition} 

\label{Elimt} For every step $s_k^n$, there exists a subsequence of $\{ ( u^n_{k,i} , v^n_{k,i} ) \}_{i\in\mathbb{N}}$ converging strongly in $W^{1,p}(\Omega,\mathbb{R}^2) \times H^1(\Omega)$; the limit, denoted $( u_k^n , v_k^n)$, is a separate minimizer for $\F$, i.e., 
\begin{equation*} 
\begin{cases}
	u_k^n \in \argmin \{  \F  (  u , v_k^n ) : u \in \U (s_k^n)
\} , 
\\
         v_k^n \in \argmin \{  \F ( u_k^n , v ) : v \in \V
\mbox{ with } v \le v^n_{k,0} = v^n_{k-1}\} . 
\end{cases} 
\end{equation*} 

\end{proposition}

\proof
By minimality
\begin{align*}
    0\leq  \F(u^n_{k,i+1},v^n_{k,i+1})\leq\F(u^n_{k,i+1},v^n_{k,i})\leq \F(u^n_{k,i},v^n_{k,i})\leq \dots \leq \F(u^n_{k,0},v^n_{k,0}):=C. 
\end{align*}
Now, since $\F ( \cdot, u)$ is quadratic we have that
\begin{align}  \label{sepq}\F(u^n_{k,i+1},v^n_{k,i})=&\ \F(u^n_{k,i+1},v^n_{k,i+1})+ \partial_v \F(u^n_{k,i+1},v^n_{k,i+1})[v^n_{k,i}-v^n_{k,i+1}]\\&+\tfrac12 \partial_{vv}^2 \F(u^n_{k,i+1},v^n_{k,i+1})[v^n_{k,i}-v^n_{k,i+1},v^n_{k,i}-v^n_{k,i+1}]. \nonumber
\end{align}

Again by minimality $\partial_v \F(u^n_{k,i+1},v^n_{k,i+1})[v^n_{k,i}-v^n_{k,i+1}]\geq 0$, while
\begin{align*}  \partial_{vv}^2 \F(u^n_{k,i+1},v^n_{k,i+1})[v^n_{k,i}-v^n_{k,i+1},v^n_{k,i}-v^n_{k,i+1}]=\int_{\Omega} 2 (v^n_{k,i}-v^n_{k,i+1})^2W_+( \strain ( u^n_{k,i+1}) )dx\\+G_c \int_{\Omega}(v^n_{k,i}-v^n_{k,i+1})^2+|\nabla v^n_{k,i}-\nabla v^n_{k,i+1}|^2dx \geq G_c \norm{v^n_{k,i}-v^n_{k,i+1}}^2_{H^1}.
\end{align*}
Therefore, from \eqref{sepq}
\begin{align*}
\tfrac12\norm{v^n_{k,i}-v^n_{k,i+1}}^2_{H^1}\leq \F(u^n_{k,i+1},v^n_{k,i})-\F(u^n_{k,i+1},v^n_{k,i+1}) \leq\F(u^n_{k,i},v^n_{k,i})   -\F(u^n_{k,i+1},v^n_{k,i+1}), 
\end{align*} and summing over $i$ it follows that
\begin{align*}
\tfrac12\sum_{i=0}^{+\infty}\norm{v^n_{k,i}-v^n_{k,i+1}}^2_{H^1}\leq \sum_{i=0}^{+\infty} \F(u^n_{k,i},v^n_{k,i})-\F(u^n_{k,i+1},v^n_{k,i+1}) \le \F(u^n_{k,0},v^n_{k,0}) \leq C. 
\end{align*} As a consequence $\|v^n_{k,i}-v^n_{k,i+1}\|^2_{H^1}\to 0$ as $i\to 0$ and from the continuous embedding of $H^1(\Omega)$ in $L^r(\Omega)$, for every $1\leq r<+\infty$, we have 
\begin{equation}\label{convv}
    \|v^n_{k,i}-v^n_{k,i+1}\|^2_{L^r}\to 0.
\end{equation} 

Since the energies $\F(u^n_{k,i},v^n_{k,i})$ are uniformly bounded, by coercivity the sequence $\{ v^n_{k,i} \}$ turns out to be bounded in $H^1(\Omega)$, and hence there exists a subsequence $\{v^n_{k,i_j}\}_{j\in\mathbb{N}}$ weakly converging in $H^1(\Omega)$ to a limit, that we call $v_{k}^n$.  By compact embedding this subsequence strongly converges in $L^r(\Omega)$ for every $1\leq r<+\infty$; as a consequence of \eqref{convv}, also $\{v^n_{k,i_j-1}\}_{j\in\mathbb{N}}$ converges to the same function. By definition
\begin{equation*}
    u^n_{k,i_j}\in \argmin\{\F(u,v^n_{k,i_j-1})\,:\, u\in \U(s^n_{k})\}. 
\end{equation*} 
By Lemma \ref{le.KRZ2.5}, the subsequence $\{u^n_{k,i_j}\}_{j\in\mathbb{N}}$ converges in $W^{1, p}(\Omega,\mathbb{R}^2)$ to $u^n_k \in \argmin\{\F(u,v^n_k)\,:\;u\in \U(s_k^n)\}$. 
Analogously, since
\begin{equation*}
    v^n_{k,i_j}\in \argmin\{\F(u^n_{k,i_j},v)\,:\, v\leq v^n_{k}\},
\end{equation*} 
by Lemma \ref{lem_dk_est_z_1WEAK}  the subsequence $\{v^n_{k,i_j}\}_{j\in\mathbb{N}}$ converges in $H^1(\Omega)$ to 
$    v^n_{k} \in \argmin\{\F(u^n_{k},v)\,:\, v\leq v^n_{k}\}$, which concludes the proof. 

\qed

%% file: appendix.tex
\section{Kuratowski convergence}
\label{Kur}

Kuratowski convergence is a notion of convergence for sequences of subsets of a topological space \cite{Kura}. 
Here, we briefly present the notion of Kuratowski limit of a sequence of subsets of a metric space $(X,d)$.
\begin{definition}Given a sequence of subsets $\{A_n\}_{n\in \mathbb{N}}$ of $X$, the Kuratowski limit inferior of $A_n$ as $n\to +\infty$ is:
\begin{align*}
    Li A_n:=\{x\in X\,:\,\limsup_{n\to+\infty} d(x,A_n)=0\}=\{x\in X\,:\,\forall \,n\in \mathbb{N}, \, \exists\,a_n \in A_n\,:\,\lim_{n\to+\infty}d(x,a_n)=0\}
\end{align*}
while the Kuratowski limit superior of $A_n$ as $n\to +\infty$ is:
\begin{align*}
    Ls A_n:=\{x\in X\,:\,\liminf_{n\to+\infty} d(x,A_n)=0\}=\{x\in X\,:\,\exists \,\{n_i\}_{i\in\mathbb{N}} \text{ and } a_{n_i} \in A_{n_i} \,:\,\lim_{i\to+\infty}d(x,a_{n_i})=0\}
\end{align*}
In general $Li A_n \subset Ls A_n$ and if they are equal, the common set is called Kuratowski limit of $A_n$.   
\end{definition} 

Note that both $Li A_n$ and $Ls A_n$ are closed subsets of $X$. Moreover, we have the following compactness result.

\begin{theorem} Let $(X,d)$ be a separable metric space and let $\{ A_n \}_{n \in \mathbb N}$ be a sequence of closed sets. There exists a subsequence of $\{ A_n \}_{n \in \mathbb N}$ converging in the sense of Kuratowski. 
\end{theorem}